\documentstyle[preprint,prb,aps]{revtex} 
\tighten 
 
\def\RR{\hbox{I \kern -.56em R}}                                
\def\CC{\hbox{\vrule height6.7pt width0.25pt \kern-.3em C}}     
\def\ZZ{\hbox{Z\kern-.33em Z}}                                  
\def\NN{\hbox{I \kern -0.55em N}}                               
\def\bull{\vrule height 1.6ex width 1.5ex depth -.1ex} 

\begin{document}
 \draft                             

\title{
 Distributional Asymptotic Expansions of Spectral Functions\\
 and of the Associated Green Kernels}

 \author{R. Estrada}
 \address{P. O. Box 276\\ Tres R\'{\i}os\\ Costa Rica}

 \author{S. A. Fulling}
\address{Department of Mathematics\\
Texas A\&M University\\
College Station, Texas 77843-3368} 

\date{March 6, 1998} 
\maketitle

\begin{abstract}

Asymptotic expansions of Green functions and  spectral densities 
associated with partial differential operators are widely applied 
in quantum field theory and elsewhere.
 The mathematical properties of these expansions can be clarified 
and more precisely determined by means of tools from distribution 
theory and summability theory.
(These are the same, insofar as recently the classic 
 Ces\`aro--Riesz theory of summability of series and integrals has 
been given a distributional interpretation.)
 When applied to the spectral analysis of Green functions
 (which are then to be expanded as series in a parameter, usually 
the time),
these methods show:  (1)  The ``local'' or ``global'' dependence of 
the expansion coefficients on the background geometry, etc., 
 is determined by the regularity of the asymptotic expansion of the 
integrand at the origin (in ``frequency space'');
 this marks the difference between a heat kernel and a Wightman
 two-point function, for instance.
(2) The behavior of the integrand at infinity determines whether 
the expansion of the Green function is genuinely asymptotic 
 in the literal, pointwise sense, or is merely valid in a 
distributional (Ces\`aro-averaged) sense;
 this is the difference between the heat kernel and the 
Schr\"odinger kernel.
 (3) The high-frequency expansion of the spectral density itself is 
local in a distributional sense (but not pointwise).
 These observations make rigorous sense out of calculations in the 
physics literature that are sometimes dismissed as merely formal.
\end{abstract}

 \vfill 
 E-mail: {\tt fulling@math.tamu.edu, 
 restrada@cariari.ucr.ac.cr} 
\pacs{} 

\section{Introduction}\label{sec1}

The aim of this article is to study several issues related to the
small-$t$ behavior of various Green functions $G(t,x,y)$ associated to
an elliptic differential operator~$H$.
These are the integral kernels of operator-valued functions of $H$,
such as the heat operator $e^{-tH}$,
the Schr\"odinger propagator $e^{-itH}$,
various wave-equation operators such as $\cos (t{\sqrt H})$,
the operator $e^{-t{\sqrt H}}$ that solves a certain elliptic 
boundary-value problem involving $H$, etc.
All these kernels are expressed 
(possibly after some redefinitions of variables) in the form
\begin{equation}G(t,x,y) =\int_0^\infty g(t\lambda) \, dE_\lambda(x,y), 
\label{(1.1)} \end{equation}
where $E_\lambda$ is the spectral decomposition of $H$,
 and $g$ is a smooth function on $(0,\infty)$.

Each such Green function raises a set of interrelated questions,
which are illumined by the simple examples listed in the Appendix:

(i)  {\em Does $G(t,x,y)$ have an asymptotic expansion as 
$t\downarrow0$?\/}
For the heat problem, (\ref{(A.1)}),
  it is well known\cite{MP,Gr} that
 \begin{mathletters}\label{(1.2)}
 \begin{equation} K(t,x,x) \sim (4\pi t)^{-d/2} \sum_{n=0}^\infty
 a_n(x,x) t^{n/2}, \label{(1.2a)}\end{equation}
 where $d$ is the dimension of the manifold ${\cal M}$ and
 $a_0(x,x) =1$.
 Similar formulas hold off-diagonal;
 for example, if ${\cal M}\subseteq \RR^d$ and the leading term in 
$H$ is the Laplacian, then 
 \begin{equation} 
 K(t,x,y) \sim (4\pi t)^{-d/2}e^{-|x-y|^2/4t} \sum_{n=0}^\infty
 a_n(x,y) t^{n/2}. \label{(1.2b)}\end{equation}\end{mathletters}
  In the case (\ref{(A.7b)}),
 the elementary heat kernel on $\RR^1$, 
  all $a_n=0$ except the first.
 In fact, this is true also of (\ref{(A.11b)}), 
 the elementary Dirichlet heat kernel on $(0,\pi)$, 
 because as $t$ goes to $0$ 
the ratio of any other term to the largest term ($e^{-(x-y)^2/4t}$)
 vanishes faster than any power of $t$.
 In particular, therefore, the expansion (\ref{(1.2)}) for fixed
 $(x,y)\in (0,\pi)\times (0,\pi)$  does not distinguish 
between the finite region $(0,\pi)$ and the infinite region 
$\RR$.
 (However, the smallness of the two nearest image terms in 
(\ref{(A.11b)})
 is not uniform near the boundary, and hence $\int_0^\pi  K(t,x,x)$
 has an asymptotic expansion 
 $(4\pi t)^{-1/2}\sum_{n=0}^\infty A_n$
  with nontrivial higher-order terms $A_n\,$.)
 This ``locality'' property will concern us again in questions (iv) 
and (v).

 The Schr\"odinger problem, (\ref{(A.2)}),
  gives rise to an expansion (\ref{(1.4)})
 that 
is formally identical to (\ref{(1.2)})
  (more precisely, obtained from it by 
the obvious analytic continuation).\cite{Sg,dW}   
 However, it is obvious from (\ref{(A.12b)}) that this expansion 
 (which again reduces to a single term in the examples 
 (\ref{(A.8)}) and (\ref{(A.12)}))
 is  {\em not literally valid}, 
 because each image term in (\ref{(A.12b)}) 
 is exactly as large in modulus 
as the ``main'' term! 

 \goodbreak
(ii)  {\em In what sense does such an expansion correspond to an 
asymptotic expansion for 
 $E_\lambda(x,y)$ as $\lambda\to +\infty$?\/}
  Formulas (\ref{(1.2)}) would follow immediately from (\ref{(1.1)})
  if 
 \begin{equation}E_\lambda(x,y) \sim \lambda^{d/2} \sum_{n=0}^\infty 
\alpha_n \lambda^{-n/2}  \label{(1.3)}\end{equation}
with $\alpha_n$ an appropriate multiple of $a_n\,$.
 The converse implication from (\ref{(1.2)}) to (\ref{(1.3)}), 
 however, is 
 generally not valid beyond the first (``Weyl'') term.
 (For example, in (\ref{(A.11a)})
  or any other discrete eigenvector expansion
 the $E_\lambda$ is a step function; its growth is described
 by $\alpha_0$ but there is an immediate contradiction with the 
form of the higher terms in (\ref{(1.3)}).)
 It has been known at least since the work of Brownell\cite{B1,B2}
  that (\ref{(1.3)}) is, nevertheless, correct
  if somehow ``averaged'' over 
sufficiently large  intervals of the variable $\lambda$.
That is, it is valid in a certain {\em distributional\/} sense.
  H\"ormander\cite{Ho1,Ho2}   
 reformulated this principle in terms of 
  literal asymptotic expansions up to some nontrivial finite order 
  for each of the {\em Riesz means\/} of $E_\lambda\,$.
 Riesz means generalize to (Stieltjes) integrals the 
 {\em Ces\`aro sums\/} 
used to create or accelerate convergence for infinite sequences and 
series (see Section~\ref{sec2}).   

  \goodbreak
(iii) {\em If an ordinary asymptotic expansion for $G$ does not 
exist, does an expansion exist in some ``averaged'' sense?\/}    
 We noted above that the Schwinger--DeWitt expansion
 \begin{equation}U(t,x,y) \sim (4\pi it)^{-d/2} e^{i|x-y|^2/4t}
 \sum_{n=0}^\infty a_n(x,y) (it)^{n/2} \label{(1.4)}\end{equation}
is not a true asymptotic expansion under the most general 
conditions.
 Nevertheless, this expansion gives correct information for the 
purposes for which it is used by (competent) physicists. 
 Clearly, the proper response in such a situation is not to reject 
the expansion as false or nonrigorous, but to define a sense (or 
more than one) in which it is true. 
At this point we cannot go into the uses made of the
Schwinger--DeWitt expansion in renormalization in quantum field 
theory (where, actually, $H$ is a hyperbolic operator instead of 
elliptic).
We can note, however, that if $U$ is to satisfy the initial 
condition in (\ref{(A.2)}),
  then as $t\downarrow0$ the main term in 
(\ref{(A.12b)}),
  which coincides with the whole of (\ref{(A.8b)}),
  must ``approach a delta function''\negthinspace,
 while the remaining terms of (\ref{(A.12b)}) must effectively vanish
 in the context of the integral 
 $\lim_{t\downarrow0}\int_0^\pi U(t,x,y) f(y)\, dy$.
 These things happen by virtue of the increasingly rapid 
oscillations of the terms, integrated against  the fixed test 
function $f(y)$.
 That is,  this instance of (\ref{(1.4)})  is literally true when 
interpreted as a relation among distributions 
 (in the variable~$y$).
 All this is, of course, well known, but our purpose here is to 
examine it in a more general context.
 We shall show that the situation for expansions like 
 (\ref{(1.4)}) is much 
like that for (\ref{(1.3)}):
 They can be rigorously established in a Riesz--Ces\`aro sense, or, 
equivalently, in the sense of distributions in the variable~$t$.
This leaves open the next question.

 \goodbreak
(iv)  {\em   If an asymptotic expansion does not exist pointwise,
does it exist distributionally in $x$ and/or $y$;
and does the spectral expansion converge in this distributional 
sense when it does not converge classically?
What is the connection between this distributional behavior and that
in~$t$?\/}
Such formulas as (\ref{(A.8a)}), (\ref{(A.10)}), 
 (\ref{(A.12a)}), (\ref{(A.14a)}) are not 
convergent, but only summable or, at most, conditionally 
convergent.
 The Riesz--Ces\`aro theory handles the summability issue, and, as 
remarked, can be rephrased in terms of distributional behavior 
in~$t$.
 However, one suspects that such integrals or sums should be 
literally convergent in the topology of distributions on ${\cal M}$ or 
${\cal M}\times{\cal M}$.
 
 This interpretation is especially appealing in the case of the 
Wightman function 
 (see (\ref{(A.4)})--(\ref{(A.5)}), (\ref{(A.10)}), (\ref{(A.14)})).
 To calculate observable quantities such as energy density in 
quantum field theory, one expects to subtract from $W(t,x,y)$ the 
leading, singular terms in the limit $y\to x$; those terms are 
``local'' or ``universal''\negthinspace, like the $a_n$ in the heat 
kernel.
The remainder will be nonlocal but finite; it contains the 
information about physical effects caused by boundary and initial 
conditions on the field.
 (See, for instance, Ref.~\onlinecite{F89}, Chapters 5 and 9.) 
 The fact that this renormalized $W(t,x,x)$ is finite does not 
guarantee that a spectral integral or sum for it will be absolutely 
convergent.  
Technically, this problem may be handled by Riesz means or some 
other definition of summability;
 but in view of the formulation of quantum field theory in terms of 
operator-valued distributions, one expects that such summability 
should be equivalent to distributional convergence on ${\cal M}$.
 It was, in fact, this problem that originally motivated the 
present work and a companion paper.~\cite{F97}.

A fully satisfactory treatment of these issues cannot be limited to
 the interior of ${\cal M}$; it should take into account the special 
phenomena that occur at the boundary.
  These questions are related to the ``heat content asymptotics'' 
  recently studied by Gilkey et al.\cite{vdB,DG}
  and McAvity.\cite{McA1,McA2}
 (A longer reference list, especially of earlier work by
Van den Berg, is given by Gilkey in 
Ref.~\onlinecite{GWin}.)

  \goodbreak
(v)  {\em  Is the expansion ``local'' or ``global'' in its dependence
on $H$?\/}
 We have already encountered this issue in connection with the 
Wightman function, but it is more easily demonstrated by what we 
call the ``cylinder kernel'' $T(t,x,y)$, defined by (\ref{(A.3)}).
 Examination of (\ref{(A.9b)}) and (\ref{(A.13b)})--(\ref{(A.13c)})
  shows that $T$ has a 
nontrivial power-series expansion in $t$, which is different for 
the two cases (${\cal M} = \RR$ and $(0,\pi)$).
 (See Ref.~\onlinecite{F97}  for more detailed discussion.)
 More generally speaking, $T(t,x,x)$ differs in an essential way 
 from $K(t,x,x)$ 
in that its asymptotic expansion as $t\downarrow0$ is not uniquely 
determined by the coefficient functions (symbol) of $H$, evaluated 
at~$x$.
 $T(t,x,x)$ can depend upon boundary conditions, existence of 
closed classical paths (geodesics or bicharacteristics), 
 and other global structure of the problem.
 In terms of an inverse spectral problem, 
the asymptotic expansion of $T$ gives more 
information about the spectrum of $H$ and about $E_\lambda(x,y)$ 
 than that of $K$ does.
(Of course, the {\em exact\/} heat kernel contains, in principle,
    all the information, as it is the Laplace transform of $E_\lambda\,$.)
We shall investigate the issue of locality for a general Green 
function (\ref{(1.1)}).

 \goodbreak
 In summary, the four basic examples introduced in the Appendix
 demonstrate all possible combinations of pointwise or 
distributional asymptotic expansions with local or global 
dependence on the symbol of the operator:
 \begin{equation}\vbox{\offinterlineskip
\halign{#\hfil&\vrule width1pt#&\strut\hfil#\hfil&\vrule#&\hfil#\hfil\cr
 &&\quad Pointwise\quad&&\quad Distributional\quad\cr
 \omit&height2pt&\omit&height2pt&\omit\cr
 \noalign{\hrule height 1pt}
 \omit&height2pt&\omit&height2pt&\omit\cr
 Local&&Heat&&Schr\"odinger\cr
 \omit&height2pt&\omit&height2pt&\omit\cr
 \noalign{\hrule}
 \omit&height2pt&\omit&height2pt&\omit\cr
 Global\quad&&Cylinder&&Wightman\cr}
 }\end{equation}

In this paper we show that the answers to questions 
 (i) and (iii),
 and the distinction between the columns of the table above, 
  are determined by
the behavior of $g$ at infinity:
If 
\begin{equation}g^{(n)}(t) = O\bigl(t^{\gamma-n}\bigr) \quad{\rm as\ }
t\to +\infty \quad{\rm for\ some\ } \gamma\in\RR  
 \label{(1.5)}\end{equation}
(i.e., $g$ has at infinity the behavior characterizing the test-function
space $\cal K$ --- see Sec.~\ref{sec2}),
 then the answer to (i) is {\em Yes}.
On the other hand, when $g$ is of slow growth at infinity but does not
necessarily belong to $\cal K$, then the expansion holds in the 
distributional sense mentioned in~(iii).

The answer to (v),
 and the distinction between the rows of the table, 
  depend on the behavior
of $g$ at the origin.
If $g(t)$ has an expansion of the form 
$\sum_{n=0}^\infty a_n t^n$ as $t\downarrow0$
(even in the distributional sense) then the expansion of $G(t,x,y)$
is local.
However, if the expansion of $g(t)$ contains fractional powers,
logarithms, or any other term, then the locality property is lost.
This subject is treated from a different point of view in 
Ref.~\onlinecite{F97}.

  We hope to return to question (iv) in later work.

Our basic tool is the study of the {\em distributional\/}
behavior of the spectral density
$e_\lambda = dE_\lambda/d\lambda$ of the operator $H$ as 
$\lambda \to \infty$.
We are able to obtain a quite general expansion of $e_\lambda$ when
$H$ is self-adjoint.
Using the results of a previous paper,\cite{E96}
  one knows that
distributional expansions are equivalent to expansions of 
Ces\`aro--Riesz means.
Thus our results become an extension of those of 
 H\"ormander.\cite{Ho1,Ho2}
 They sharpen and complement previous publications by one of 
 us.\cite{F1,F3,F97}

The other major tool we use is an extension of the ``moment
asymptotic expansion'' to distributions, as explained in 
 Section~\ref{sec5}.

The plan of the paper is as follows. 
 In Section~\ref{sec2} we give some results
from Ref.~\onlinecite{E96}
  that play a major role in our
analysis. 
 In particular we introduce the space of test functions
${\cal K}$ and its dual ${\cal K}'$, the space of distributionally small
generalized functions.

In the third section we consider the distributional asymptotic expansion
of spectral decompositions and of spectral densities. Many of our
results hold for general self-adjoint operators on a Hilbert space, and
we give them in that context. We then specialize to the case of a
pseudodifferential operator acting on a manifold and by exploiting
the pseudolocality of such operators we are able to show that the
asymptotic behavior of the spectral density of a pseudodifferential
operator has a local character {\it in the Ces\`aro sense}. That such
spectral densities have a local character in ``some sense'' has been 
known for years;\cite{FNW,F1,F2,F3} here we provide a precise meaning to
this locality property.

In the next section we consider two model examples for the
asymptotic expansion of spectral densities. Because of the local behavior,
they are more than examples, since they give the asymptotic development
of any operator locally equal to one of them.

In Section~\ref{sec5} we show that the moment asymptotic expansion,
  which is the 
basic building block in the asymptotic expansion of series and
integrals,\cite{EK94}
  can be generalized to distributions, giving expansions
that hold in an ``averaged'' or distributional sense explaining,
for instance, the small-$t$ behavior of the Schr\"odinger propagator.

In the last two sections we apply our machinery to the study of the 
asymptotic expansion of general Green kernels. 
 In Section~\ref{sec6} we show
that the small-$t$ expansion of a propagator $g(tH)$ that corresponds
to a function $g$ that has a Taylor-type expansion at the origin is
local and that it is an ordinary or an averaged expansion depending
on the behavior of $g$ at infinity: 
 If $g\in{\cal K}$ then the regular
moment asymptotic expansion applies, while if $g\not\in{\cal K}$ then
the ``averaged'' results of Section~\ref{sec5} apply. 
 In the last section we 
consider the case when $g$ does not have a Taylor expansion at
the origin and show that in that case $g(tH)$ has a global expansion,
which depends on such information as boundary conditions.

Some applications of both of the main themes of this paper have been
made elsewhere,\cite{EGV} most notably a mathematical sharpening of
the work of Chamseddine and Connes\cite{CC} on a ``universal bosonic
functional''\negthinspace.

 \section{Preliminaries}\label{sec2}

The principal tool for our study of the behavior of spectral functions
and of the associated Green kernels is the distributional theory of
asymptotic expansions, as developed by several 
 authors.\cite{VDZ,P,EK90,EK94}
The main idea is that one may obtain the ``average" behavior of a
function, in the Riesz or Ces\`aro sense, by studying its parametric
or distributional behavior.\cite{E96}

In this section we give a summary of these results. We also set
the notation for the spaces of distributions and test functions used.

If ${\cal M}$ is a smooth manifold, 
 then ${\cal D}({\cal M} )$ is the space of 
compactly supported smooth functions on ${\cal M}$, equipped with the 
standard Schwartz topology.\cite{EK94,S,K}
Its dual, ${\cal D}'({\cal M} )$, 
is the space of standard distributions on ${\cal M}$. 
 The space ${\cal E}({\cal M} )$ is 
the space of all smooth functions on ${\cal M}$, 
 endowed with the topology of 
uniform convergence of all derivatives on compacts. 
 Its dual, ${\cal E}'({\cal M} )$,
can be identified with the subspace of ${\cal D}'({\cal M} )$
formed by the compactly supported distributions. Naturally the 
two constructions coincide if ${\cal M}$ is compact.

The space ${\cal S}'(\RR^n)$ consists of the tempered distributions 
 on $\RR^n$. 
 It is the dual of the space of rapidly decreasing smooth 
 functions ${\cal S}(\RR^n)$;
a smooth function $\phi$ belongs to ${\cal S}(\RR^n)$ if $D^\alpha\phi(x) 
=o(|x|^{-\infty})$ as $|x|\to\infty$, for each $\alpha\in\NN^n$. 
Here we use the usual notation,
$D^\alpha =\partial^{|\alpha |}/\partial x_1^{\alpha_1}\cdots \partial
x_n^{\alpha_n}$, $|\alpha |=\alpha_1+\cdots + \alpha_n$; $o(x^{-\infty})$
means a quantity that is $o(x^{-\beta})$ for all $\beta\in\RR$.

A not so well known pair of spaces that plays a fundamental role in our
analysis is ${\cal K}(\RR^n)$ and ${\cal K}'(\RR^n)$. 
 The space ${\cal K}$ was introduced in Ref.~\onlinecite{GLS}.  
 A smooth function $\phi$ belongs to ${\cal K}_q$ if $D^\alpha\phi(x)=
O(|x|^{q-|\alpha |})$ as $|x|\to\infty$ for each $\alpha\in\NN^n$. 
 The
space ${\cal K}$ is the inductive limit of the spaces ${\cal K}_q$ 
as $q\to\infty$. 

Any distribution $f\in{\cal K}'(\RR)$ satisfies the 
 {\em moment asymptotic expansion},
\begin{equation}
f(\lambda x)\sim\sum_{j=0}^\infty {(-1)^j\mu_j\delta^{(j)}(x)\over
j!\,\lambda^{j+1}}\quad  {\rm as}\ \lambda\to\infty,
\label{(2.1)}
\end{equation}
 where $\mu_j=\langle f(x),x^j\rangle$ are the moments of $f$. 
 The interpretation of (\ref{(2.1)})
  is in the topology of the space ${\cal K}'$;
observe, however, that there is an equivalence between weak and strong
convergence of one-parameter limits in spaces of distributions,
such as ${\cal K}'$.

The moment asymptotic expansion does not hold for general distributions
of the spaces ${\cal D}'$ or ${\cal S}'$. 
 Actually, it was shown recently\cite{E96}
that any distribution $f\in{\cal D}'$ that satisfies the moment expansion 
(\ref{(2.1)})
  for some sequence of constants $\{\mu_j\}$ must belong to ${\cal K}'$
(and then the $\mu_j$ are the moments).

There is still another characterization of the elements of ${\cal K}'$.
  They are
precisely the distributions of rapid decay at infinity {\em in the
Ces\`aro sense}.
  That is why the elements of ${\cal K}'$ are referred to as
{\em distributionally small}.

The notions of Ces\`aro summability of series and integrals 
 are well known.\cite{Ha2}
In Ref.~\onlinecite{E96}  this theory 
 is generalized to general distributions. 
 The
generalization includes the classical notions as particular cases, 
 since the
behavior of a sequence $\{a_n\}$ as $n\to\infty$ can be studied by 
studying the generalized function 
 $\sum_{n=0}^\infty\, a_n\,\delta (x-n)$.
The basic concept is that of the order symbols in the Ces\`aro sense:
  Let
$f\in{\cal D}'(\RR)$ and let $\beta\in\RR\setminus \{-1,-2,-
3,\ldots\}$; 
 we say that
\begin{equation}
f(x)=O(x^\beta )\quad (C)\quad {\rm as}\ x\to\infty,
\label{(2.2)}\end{equation}
if there exists $N\in\NN$, 
 a function $F$ whose $N$th derivative is $f$,  and a
polynomial $p$ of degree $N-1$ such that $F$ is locally integrable for
$x$ large and the ordinary relation
\begin{equation}
F(x)=p(x)+O(x^{\beta +N})\quad {\rm as}\ x\to\infty
\label{(2.3)} 
\end{equation}
holds. 
 The relation $f(x)=o(x^\beta )\ \ (C)$ is defined similarly by replacing
the big $O$ by the little $o$ in (\ref{(2.3)}).

Limits and evaluations can be handled by using the order relations.
In particular, $\lim_{x\to\infty}f(x)=L\ \ (C)$ means that
$f(x)=L+o(1)\ \ (C)$ as $x\to\infty$. 
 If $f\in{\cal D}'$ has support bounded on the
left and $\phi\in {\cal E}$, then in general the evaluation
$\langle f(x),\phi (x)\rangle$ does not exist, but we say that it has the
value $S$ in the Ces\`aro sense if $\lim_{x\to\infty}G(x)=S\ \ (C)$,
where $G$ is the primitive of $f\phi$ with support bounded on the left.
The Ces\`aro interpretation of evaluations $\langle f(x),\phi(x)\rangle$
with $\mathop{\rm supp}{f}$ bounded on the right is similar, 
 while the general case
can be considered by writing $f=f_1+f_2$, 
 with $\mathop{\rm supp}{f_1}$ bounded on
the left and $\mathop{\rm supp}{f_2}$ bounded on the right.

The main result that allows one to obtain the Ces\`aro behavior from the
parametric behavior is the following.

\smallskip
{\bf Theorem 2.1.} 
{\em  Let $f$ be in ${\cal D}'$ with support bounded on the left. If
$\alpha >-1$, then
\begin{equation}
f(x)=O(x^\alpha )\quad (C)\quad {\rm as}\ x\to\infty
\label{(2.4)}\end{equation}
if and only if
\begin{equation}
f(\lambda x)=O(\lambda^\alpha )\quad {\rm as}\ \lambda\to\infty
\label{(2.5)}\end{equation}
distributionally.

When $-(k+1)>\alpha >-(k+2)$ for some $k\in\NN$, 
 (\ref{(2.4)}) holds if and
only if there are constants $\mu_0,\ldots ,\mu_k$ such that
\begin{equation}
f(\lambda x)=\sum_{j=0}^k{(-1)^j\mu_j\delta^{(j)}(x)\over j!\,\lambda^{j+1}}
+O(\lambda^\alpha )
\label{(2.6)}\end{equation}
distributionally as $\lambda\to\infty$.
 }\smallskip

{\em Proof:\/} See Ref.~\onlinecite{E96}.  $\;\;\bull$
 \smallskip

The fact that the distributions that satisfy the moment asymptotic 
expansion are exactly those that satisfy $f(x)=O(x^{-\infty})\ \ (C)$ 
follows from the theorem by letting $\alpha\to -\infty$. 
 Thus the 
elements of ${\cal K}'$ are the distributions of rapid 
distributional decay at infinity in the Ces\`aro sense. 
 Hence the 
space  ${\cal K}'$ is a distributional analogue of ${\cal S}$. 
 We apply this idea in Section~\ref{sec5}, where we build a duality 
between ${\cal S}'$ and ${\cal K}'$. 

Another important corollary of the theorem is the fact that one can 
relate the $(C)$ expansion of a generalized function and its 
parametric expansion in a simple fashion.
  Namely, if
$\{\alpha_j\}$ is a sequence with $\Re e\;\alpha_j\searrow -\infty$,
then  
\begin{equation}
f(x)\sim \sum_{j=0}^\infty a_jx^{\alpha_j}\quad (C)\quad
  {\rm as}\ x\to\infty
\label{(2.7)}\end{equation}
if and only if
\begin{equation}
f(\lambda x)\sim\sum_{j=0}^\infty a_jg_{\alpha_j}(\lambda x)+
\sum_{j=0}^\infty {(-1)^j\mu_j\delta^{(j)}(x)\over j!\,\lambda^{j+1}}
\label{(2.8)}\end{equation}
as $\lambda\to\infty$, where the $\mu_j$ are the (generalized)
moments of $f$ and where 
\begin{equation}
g_\alpha (x) =x_+^\alpha\quad{\rm if}\ \alpha\not= -1,-2,-3,\ldots,
\label{(2.9)}\end{equation}
while in the exceptional cases $g_\alpha$ is a finite-part 
 distribution:\cite{EK94}
\begin{equation}
g_{-k}(x)={\cal P}.f.\,(\chi(x)x^{-k})\quad {\rm if}\ k=1,2,3,\ldots,
\label{(2.10)}\end{equation}
 $\chi$ being the Heaviside function, the characteristic function of
the interval $(0,\infty )$.
Notice that 
\begin{equation}
g_\alpha (\lambda x)=\lambda^\alpha g_\alpha (x),\quad
   \alpha\not=-1,-2,-3,\ldots,
\label{(2.11)}\end{equation}
\begin{equation}
g_{-k}(\lambda x)=
 {g_{-k}(x)\over \lambda^k}+ {(-1)^{k-1}\ln{\lambda}\,\delta^{(k-1)}
(x)\over (k-1)!\,\lambda^k}\,,\quad k=1,2,3,\ldots.
\label{(2.12)}\end{equation}
   
\section{The asymptotic expansion of spectral decompositions} 
\label{sec3} 
  
Let ${\cal H}$ be a Hilbert space and 
 let $H$ be a self-adjoint operator on ${\cal H}$, 
with domain ${\cal X}$.   
 Then $H$ admits a spectral decomposition $\{E_{\lambda} 
\}_{\lambda=-\infty}^{\infty}\,$.   
 The $\{E_{\lambda}\}$ is an increasing  
family of projectors that satisfy 
\begin{equation} I = \int_{-\infty}^{\infty} d E_{\lambda}   \, , 
\label{(3.1)}  
\end{equation} 
where $I$ is the identity operator, and 
\begin{equation} H = \int_{-\infty}^{\infty} \lambda\, d E_{\lambda}  
\label{(3.2)}  
\end{equation} 
in the weak sense, that is, 
\begin{equation} (Hx|y) = \int_{-\infty}^{\infty}  
 \lambda\, d (E_{\lambda}x|y)    ,
\label{(3.3)}  
\end{equation} 
for $x \in {\cal X}$ and $y \in {\cal H}$, 
 where $(x|y)$ is the inner product in ${\cal H}$. 
 
Perhaps more natural than the spectral function $E_{\lambda}$ is the 
spectral density $e_{\lambda}= d E_{\lambda}/d \lambda$. 
  This spectral density 
does not have a pointwise value for all $\lambda. $  Rather, it should be  
understood as an operator-valued distribution, 
 an element of the space ${\cal D}'(\RR, L({\cal X},{\cal H})). $  
 Thus (\ref{(3.1)})--(\ref{(3.2)}) become 
\begin{equation} I =  \langle e_{\lambda}, 1 \rangle  
\label{(3.4)}  
\end{equation} 
\begin{equation} H =  \langle e_{\lambda}, \lambda \rangle    , 
\label{(3.5)}  
\end{equation} 
where $\langle f(\lambda), \phi(\lambda) \rangle$  
 is the evaluation of a distribution 
$f(\lambda)$ on a test function $\phi(\lambda)$. 
 
The spectral density $e_{\lambda}$ can be used to build a functional 
calculus for the operator $H. $  
 Indeed, if $g$ is continuous and with compact support  
in $\RR$ then we can define the operator $g(H) \in L({\cal X},{\cal H})$  
 (extendible to $L({\cal H},{\cal H})$) by 
\begin{equation} g(H) =  \langle e_{\lambda}, g(\lambda) \rangle    . 
\label{(3.6)}  
\end{equation} 
One does not need to assume $g$ of compact support in (\ref{(3.6)}), 
  but in  a contrary case
the domain of $g(H)$ is not ${\cal X}$ but the subspace ${\cal N}_g$ 
 consisting of the $x  
\in {\cal H}$ for which the improper integral 
 $\langle (e_{\lambda} x|y), g(\lambda)  
\rangle$ converges for all $y \in {\cal H}$. 
  
One can even define $f(H)$ when $f$ is a distribution such that the 
evaluation $\langle e_{\lambda}, f(\lambda) \rangle$ is defined.   
 For instance, if $E_{\lambda}$ 
is continuous at $\lambda=\lambda_0$ 
 then $E_{\lambda_0}=\chi(\lambda_{0}-H)$ 
where $\chi$ is again the Heaviside  function. 
  Differentiation yields 
the useful symbolic formula 
\begin{equation} e_{\lambda} = \delta(\lambda-H)    . 
\label{(3.7)}  
\end{equation} 
 
Let ${\cal X}_n$ be the domain of $H^n$ and let  
 ${\cal X}_{\infty}=\displaystyle\bigcap_{n=1}^{\infty} {\cal X}_n\,. $ 
Then 
\begin{equation} \langle e_{\lambda}, \lambda^n \rangle =H^n 
\label{(3.8)}  
\end{equation} 
in the space $L({\cal X}_{\infty},{\cal H}). $ 
  But, as shown recently,\cite{E96}
a distribution $f \in {\cal D}'(\RR)$ whose moments  
 $\langle f(x),x^n \rangle$, $n \in \NN$, 
all exist belongs to ${\cal K}'(\RR)$, that is, is distributionally small.
   Hence,  
$e_{\lambda}\,$, as a function of $\lambda$, belongs to the space  
 ${\cal K}'(\RR,L({\cal X}_{\infty},  
{\cal H})). $   
 Therefore, the asymptotic behavior of $e_{\lambda \sigma}\,$, as $\sigma   
\rightarrow \infty$, can be obtained by using the moment asymptotic expansion:  
\begin{equation} e_{\lambda \sigma} \sim  
 \sum_{n=0}^{\infty} \frac{(-1)^n H^n \delta^{(n)}(\lambda   
\sigma)}{n!} \quad \mbox{as}\ \sigma \rightarrow \infty    ,  
\label{(3.9)}   
 \end{equation}  
while $e_{\lambda}$ vanishes to infinite order at infinity in the  
 Ces\`aro sense:  
\begin{equation} e_{\lambda} = 
 o(|\lambda|^{-\infty}) \quad (C)   \quad \mbox{as}\
 |\lambda| \rightarrow \infty    .  
\label{(3.10)}   
\end{equation}  
  
The asymptotic behavior of the spectral function $E_{\lambda}$ is obtained by  
integration of (\ref{(3.9)}) and by recalling that 
$\displaystyle\lim_{\lambda  \rightarrow -\infty}  
E_{\lambda}=0  , \, 
 \displaystyle\lim_{\lambda  \rightarrow \infty}E_{\lambda}=I. $  
We obtain  
\begin{equation} E_{\lambda} \sim \chi(\lambda \sigma)I + \sum_{n=0}^{\infty} 
\frac{(-1)^{n+1}   
H^{n+1} \delta^{(n)}(\lambda \sigma)}{(n+1)!} 
 \quad \mbox{as}\    
\sigma \rightarrow \infty.
\label{(3.11)}   
\end{equation}  
   Similarly, the Ces\`aro behavior is   
given by  
\begin{equation} E_{\lambda} = I + o(\lambda^{-\infty}) \quad
(C)\quad
 \mbox{as}  \
 \lambda \rightarrow \infty    ,  
\label{(3.12)}   
\end{equation}  
\begin{equation} E_{\lambda} = o(|\lambda|^{-\infty}) \quad
(C)\quad
 \mbox{as}  
\ \lambda \rightarrow -\infty    .  
\label{(3.13)}   
\end{equation}  
  
 These formulas are most useful when $H$ is an   
unbounded operator. 
  Indeed, if $H$ is bounded, with domain ${\cal X}={\cal H}$, then   
$e_{\lambda}=0$ for $\lambda>\| H\|$ and  
 $E_{\lambda}=0$ for $ \lambda<-\| H\|, \; 
E_{\lambda}=I$ for $ \lambda>\| H\|, $  
so (\ref{(3.10)}), (\ref{(3.12)}), and (\ref{(3.13)}) are trivial  
in that case.  
  
In the present study we are mostly interested in the case when $H$  
 is an elliptic  
differential operator with smooth coefficients defined on a  
 smooth manifold ${\cal M}$.   
Usually ${\cal H}=L^2({\cal M})$ and ${\cal X}$ is the domain 
corresponding to  
 the introduction of  
suitable boundary conditions.  
 Usually the operator $H$ will be positive, but  
at present we shall just assume $H$ to be self-adjoint. 
  
In this case the space ${\cal D}({\cal M})$ of test functions on  
 ${\cal M}$ is a subspace of ${\cal X}_{\infty}\,$.   
Observe also that the operators $K$ acting on ${\cal D}({\cal M})$ 
  can be realized as distributional  
kernels $k(x,y)$ of ${\cal D}'({\cal M} \times {\cal M})$ by  
\begin{equation} (K \phi)(x) = \langle k(x,y), \phi(y) \rangle_y \, .  
\label{(3.14)}   
\end{equation}  
In particular, 
 $\delta(x-y)$ is the kernel corresponding to the identity $I$, and 
$H \delta(x-y)$ is the kernel of $H. $  
The spectral density $e_{\lambda}$ also  
has an associated kernel $e(x,y;\lambda)$, an element of
 ${\cal D}'(\RR, {\cal D}'({\cal M} \times {\cal M} )). $  
 Since $H$ is elliptic it follows 
 that $e(x,y;\lambda)$ is smooth in $(x,y)$.  
 {\em Warning:\/} Much of the literature uses ``$e(x,y;\lambda)$''
 for what we call ``$E(x,y;\lambda)$''\negthinspace.
  
The expansions (\ref{(3.9)})--(\ref{(3.13)})
  will hold in ${\cal X}_{\infty}$ and thus, 
 a fortiori, in   
${\cal D}({\cal M} ). $  Hence  
\begin{equation} e(x,y;\lambda \sigma) \sim 
 \sum_{n=0}^{\infty} \frac{(-1)^n H^n \delta(x-y)  
\delta^{(n)}(\lambda \sigma)}{n!}  \quad \mbox{as} \ \sigma  
\rightarrow \infty    ,   
\label{(3.15)}   
\end{equation}  
\begin{equation} E(x,y;\lambda \sigma) \sim \chi(\lambda \sigma) 
\delta(x-y)   + \sum_{n=0}^{\infty}   
\frac{(-1)^{n+1} H^{n+1} \delta(x-y) \,\delta^{(n)}(\lambda \sigma)} 
 {(n+1)!} \quad \mbox{as} \ \sigma  
\rightarrow \infty    ,   
\label{(3.16)}   
\end{equation}  
in the space ${\cal D}'_{\lambda}(\RR,{\cal D}'_{xy}({\cal M}
  \times {\cal M} ))$.  
 Furthermore, 
\begin{equation} e(x,y;\lambda) = o(|\lambda|^{-\infty}) \quad 
(C)\quad
 \mbox{as}  
\ |\lambda| \rightarrow \infty ,  
\label{(3.17)}   
\end{equation}  
 \begin{mathletters}\label{(3.18)}
\begin{equation} E(x,y;\lambda) = \delta(x-y) + o(\lambda^{-\infty}) 
 \quad (C) \quad
 \mbox{as}  
\ \lambda \rightarrow \infty    ,  
\label{(3.18a)}   
\end{equation}  
\begin{equation} E(x,y;\lambda) = o(|\lambda|^{-\infty}) \quad  (C) 
\quad
  \mbox{as}  
\ \lambda \rightarrow -\infty    ,  
\label{(3.18b)}   
\end{equation}  \end{mathletters}
in the space ${\cal D}'({\cal M} \times {\cal M} )$.  
  
Actually, an easy argument shows that the expansions also hold 
distributionally  
in one variable and pointwise in the other. 
  (For instance, (\ref{(3.17)}) says that if $y$  
is fixed and $\phi \in {\cal D}({\cal M} )$ then  
 $\langle e(x,y;\lambda), \phi(x) \rangle  
=o(|\lambda|^{-\infty}) \ \ (C) \ \ \mbox{as} \; |\lambda| 
\rightarrow \infty$.)  
  
That the expansions cannot hold pointwise in both variables $x$ and 
$y$ should  be clear since we cannot set $x=y$ in the distribution 
$\delta(x-y). $   
 And indeed, 
$e(x,x;\lambda)$ is {\em not\/} distributionally small.  
However, as we now show,  
the expansions are valid pointwise outside of the diagonal of  
 ${\cal M} \times{\cal M}$.  
  
Indeed, let $U, V$ be open sets with $U \cap V= \emptyset. $   
 If $f \in {\cal D}'({\cal M})$  
and $\phi \in {\cal D}(\RR)$, then $\phi(H)$ is a smoothing 
  pseudodifferential operator,  
so $\phi(H)f$ is smooth in ${\cal M}$.
   Thus,  
$\langle e(x,y;\lambda),f(x)g(y)\phi(\lambda) \rangle =
  \langle \phi(H)f(x),  g(x)\rangle$ 
 is well-defined if 
 $f \in {\cal D}'({\cal M} ), \mathop{\rm supp} f \subseteq U, 
  g \in {\cal D}'({\cal M}), 
\mathop{\rm supp} g \subseteq V. $   
 Therefore $e(x,y;\lambda)$ belongs to  
 ${\cal D}'(\RR,{\cal E}(U \times V))$.   
But  
\begin{equation} \langle e(x,y;\lambda),f(x)g(y) \lambda^n \rangle
  = \langle H^n f(x),  
g(x)\rangle = 0    ,  
\label{(3.19)}   
\end{equation}  
thus $e(x,y;\lambda)$ actually belongs to 
 ${\cal K}'(\RR,{\cal E}(U \times V))$; 
  that is, it is a  
distributionally small distribution in that space whose moments vanish. 
   Therefore  
\begin{equation} e(x,y;\lambda \sigma) = o(\sigma^{-\infty})  
  \quad \mbox{as}  
\ \sigma \rightarrow \infty    ,  
\label{(3.20)}   
\end{equation}  
\begin{equation} E(x,y;\lambda \sigma) = \chi(\lambda \sigma) 
\delta(x-y) + o(\sigma^{-\infty})  
\quad \mbox{as} \ \sigma \rightarrow \infty    ,  
\label{(3.21)}   
\end{equation}  
in the space ${\cal K}'(\RR,{\cal E}(U \times V)). $   
 Similarly, (\ref{(3.17)})--(\ref{(3.18)}) also hold in ${\cal E}(U  
\times V). $   
 Convergence in ${\cal E}(U \times V)$ implies  
 pointwise convergence on $U\times V$, but 
it gives more; namely, it gives uniform convergence of all  
 derivatives on compacts.  
Thus (\ref{(3.17)}), (\ref{(3.18)}), (\ref{(3.20)}), and 
(\ref{(3.21)}) 
 hold uniformly on compacts of $U \times  
V$ and the expansion can be differentiated as many times as we 
please with  respect to $x$ or $y$.  

 \smallskip  
{\bf Example.} Let $H y = -y''$ considered on the domain  
 ${\cal X}=\{ y \in C^2[0,\pi]:  
y(0)=y(\pi)=0 \}$ in $L^2[0,\pi]. $  
 The eigenvalues are $\lambda_n=n^2$, 
$n=1,2,3,\ldots$, with normalized eigenfunctions 
$\phi_n(x)=\sqrt{\frac{2}{\pi}}  
\sin nx. $ 
  Therefore,  
\begin{equation} e(x,y;\lambda) =  
 \frac{2}{\pi} \sum_{n=1}^{\infty} \sin nx \sin ny \, \delta(\lambda-  
n^2)    ,  
\label{(3.22)}   
\end{equation}  
where $0<x<\pi$, $0<y<\pi$. 
   Then  
\begin{equation} \frac{2}{\pi} \sum_{n=1}^{\infty} \sin nx \sin ny 
\, \delta(\lambda \sigma-n^2)  
\sim \sum_{j=0}^{\infty} 
 \frac{\delta^{(2j)}(x-y)\delta^{(j)}(\lambda)}{j! \,
\sigma^{j+1}} \quad \mbox{as} \ \sigma \rightarrow \infty   
\label{(3.23)}   
\end{equation}  
in ${\cal D}'(\RR,{\cal D}'((0,\pi) \times (0,\pi)) )$, while  
\begin{equation}  \frac{2}{\pi} \sum_{n=1}^{\infty} \sin nx \sin ny 
\, \delta(\lambda \sigma-n^2)  
= o(\sigma^{-\infty}) \quad
\mbox{as} \ \sigma \rightarrow \infty   
\label{(3.24)}   
\end{equation}  
if $x$ and $y$ are fixed, $x \neq y. $  On the other hand,  
\begin{equation} e(x,x;\lambda) =
  \frac{2}{\pi} \sum_{n=1}^{\infty} \sin^2 nx \,
  \delta(\lambda-  n^2)    ,  
\label{(3.25)}   
\end{equation}  
thus if $0<x<\pi$,  
\begin{eqnarray*}  
e(x,x;\lambda \sigma) & = & \frac{1}{\pi} \sum_{n=1}^{\infty} (1-
\cos 2nx) \,
\delta(\lambda \sigma-n^2) \\  
& = & \frac{1}{\pi} \sum_{n=1}^{\infty} \delta(\lambda \sigma-n^2) +   
\frac{1}{2\pi \sigma} \,\delta(\lambda) + o(\sigma^{-\infty})
\quad\mbox{as} \
  \sigma \rightarrow \infty,  
\end{eqnarray*}  
because the generalized function  
 $\sum_{n=1}^\infty \cos{2nx}\,\delta (\lambda -n^2)$  
is distributionally small if $0<x<\pi$, with moments 
$\mu_0=-1/2$ and $\mu_k=0$  
for $k\geq 1$, since\cite{E95}  
\begin{eqnarray*}  
\sum_{n=1}^\infty \cos{2nx} &=& -\,{1\over 2}\quad  (C),  
\\  
\sum_{n=1}^\infty n^{2k}\cos{2nx} &=& 0 \quad (C),
 \qquad k=1,2,3,\ldots.  
\end{eqnarray*}  
But (Ref.~\onlinecite{EK94}, Chapter 5)   
\begin{equation} \sum_{n=1}^{\infty} \phi(\varepsilon n^2) 
 = \frac{1}{2 \varepsilon^{1/2}}  
\int_{0}^{\infty} u^{-1/2} \phi(u)\,du - \frac{1}{2} \phi(0) 
+ o(\varepsilon^{\infty})  
\label{(3.26)}   
\end{equation}  
as $\varepsilon \rightarrow 0^{+}$ if $\phi \in {\cal S}$, thus  
\begin{equation}    
e(x,x;\lambda\sigma )={1\over 2\pi\sigma^{1/2}}\,\lambda_+^{-1/2}
+o(\sigma^{-\infty})
\quad  {\rm as}\ \sigma\to\infty.  
\label{(3.27)}  
\end{equation}  
It is then clear that  
$e(x,x;\lambda)$ is not distributionally small; rather, 
\begin{equation} e(x,x;\lambda)= \frac{1}{2\pi \lambda^{-1/2}} 
 + o(\lambda^{-\infty}) 
 \quad  
(C) \quad \mbox{as} \ \lambda \rightarrow \infty    ,  
\label{(3.28)}   
\end{equation}  
that is, $e(x,x;\lambda) \sim (1/2\pi) \lambda^{-1/2}, \; \mbox{as} \;\;   
\lambda \rightarrow \infty$, in the Ces\`aro sense. $\;\; \bull$  
\smallskip

Neither is the spectral density $e(x,y;\lambda)$  distributionally small 
at the  
boundaries, as follows from the heat content asymptotics of  
Refs.~\onlinecite{vdB,DG}.  
That there is a sharp change of behavior at the boundary can be seen  
from the behavior of the spectral density $e(x,x;\lambda)$ given by 
 (\ref{(3.25)}).    
Indeed, if $0<x<\pi$ then $e(x,x;\lambda)=(1/2\pi) \lambda^{-1/2}+  
o(\lambda^{-\infty}) \ \ (C)$, 
 but when $x=0$ or $x=\pi$ then $e(0,0;\lambda)=  
e(\pi,\pi;\lambda)=0$.  
  
It is important to observe that in the Ces\`aro or distributional 
sense, the  behavior at infinity of the spectral density 
$e(x,y;\lambda)$ depends only on  the {\em local\/} behavior of the 
coefficients of $H$. 
 That is, if $H_1$ and  
$H_2$ are two operators that coincide on the open subset $U$ of 
${\cal M}$ and if   
$e_1(x,y;\lambda)$ and $e_2(x,y;\lambda)$  
 are the corresponding spectral densities,  
then  
\begin{equation} e_1(x,y;\sigma\lambda) = 
 e_2(x,y;\sigma\lambda) + o(\sigma^{-\infty}) 
  \quad  
\mbox{as} \ \sigma \rightarrow \infty   
\label{(3.29)}   
\end{equation}  
in ${\cal D}'(U \times U)$. 
 This  follows immediately from (\ref{(3.15)}).
   In fact, 
it follows from Theorem~7.2 that  
\begin{equation} e_1(x,y;\lambda) =
  e_2(x,y;\lambda) + o(\lambda^{-\infty})  
 \quad (C) \quad
\mbox{as} \;\; \lambda \rightarrow \infty    ,  
\label{(3.30)}   
\end{equation}  
pointwise on $(x,y) \in U \times U$ (even on the diagonal!).  
  More than that, (\ref{(3.30)})  
holds in the space ${\cal E}(U \times U)$, 
so that it is uniform on compacts of $U $.  
These results are useful in connection with the 
 suggestion\cite{FNW,F1,F2,F3} 
to replace a general  
 second-order operator $H$ by another, $H_0\,,$  
 that agrees locally with $H$ and  
for which the spectral density  can be calculated.   
In the next section we treat two special classes of operators  
where this idea has been implemented.   

 \smallskip  
{\bf Example.}
  The spectral density for the operator $-y''$ on the whole real   
line is   
\begin{equation}e_1(x,y;\lambda)= {\chi(\lambda)\cos{\lambda^{1/2}(x-y)} 
 \over 2\pi\lambda^{1/2}}\,,  
\label{(3.31)}
\end{equation}  
as can be seen from (\ref{(A.7a)}) and (\ref{(A.8a)}).  
 Therefore, comparison with (\ref{(3.22)}) yields  
\begin{equation}  
{2\over\pi }\sum_{n=1}^\infty \sin{nx}\sin{ny}\,\delta (\lambda -n^2)=  
{\cos{\lambda^{1/2}(x-y)}\over 2\pi\lambda^{1/2}}+o(\lambda^{-\infty})
\quad (C)  
\quad {\rm as}\  \lambda\to\infty.  
\label{(3.32)}  
\end{equation}  
In particular, if we set $x=y$ we recover (\ref{(3.28)}).  
  
Formula (\ref{(3.32)}) 
 is uniform in compacts of $(0,\pi )\times (0,\pi )$ but ceases  
to hold as $x$ or $y$ approaches $0$ or $\pi$.  
 For instance, if $y=0$, the  
left side vanishes while [cf.~(\ref{(3.15)})]  
\begin{displaymath}  
{\chi(\sigma\lambda)\cos{(\sigma\lambda)^{1/2}x}\over 2\pi  
 (\sigma\lambda)^{1/2}}  
\sim {\delta (x)\delta (\lambda)\over \sigma}+
{\delta''(x)\delta'(\lambda)\over  
\sigma^2}+\cdots       
\end{displaymath}  
as $\sigma\to\infty$. $\;\;\bull$  
\smallskip  
  
\section{Special cases}  \label{sec4}

In this section we give two model cases for the  
 asymptotic expansion of spectral  
densities. 
  They are not just examples, since according to the results of the  
previous section, the spectral density of any operator locally 
equal to such   a model case will have the same behavior at 
infinity in the Ces\`aro sense.  
  
Let us start with a constant-coefficient elliptic operator $H$ 
defined on the  whole space $\RR^n $. 
  Then $H$ admits a unique self-adjoint extension (which  
we also denote as $H$), given as follows. 
  Let $p=\sigma(H)$ be the symbol of   
$H $ (i.e., $H = p(-i\partial)$).  
   Then the spectral function is given by 
\begin{equation} E(x,y;\lambda) = 
 \frac{1}{(2 \pi)^n} \int_{p(\xi)<\lambda} e^{i(x-y) \cdot  
\xi}\,d\xi    ,  
\label{(4.1)}   
\end{equation}  
so that the spectral density can be written as  
\begin{equation} e(x,y;\lambda) =
  \frac{1}{(2 \pi)^n} \left\langle e^{i(x-y) \cdot \xi} ,  \,
 \delta  
(p(\xi)-\lambda) \right\rangle    .  
\label{(4.2)}   
\end{equation}  
For the definition of $\delta(f(x))$ 
 see Refs.~\onlinecite{GS,dJ}.  
  
To obtain the behavior of $e(x,y;\lambda)$ as $\lambda \rightarrow 
\infty$ in  the Ces\`aro or in the distributional sense, we should 
consider the parametric  behavior of $e(x,y;\sigma \lambda)$ as 
$\sigma \rightarrow \infty$.   
 Setting  
$\varepsilon=1/\sigma$ and evaluating at a test function 
$\phi(\lambda)$, one  is led to the function  
 \begin{equation} \Phi(\varepsilon)
  = \left\langle e(x,y;\lambda),\phi(\varepsilon \lambda)  
\right\rangle_{\lambda}    \, .  
\label{(4.3)}   
\end{equation}  
But in view of (\ref{(4.2)}) we obtain  
\begin{equation} \Phi(\varepsilon) =
  \frac{1}{(2\pi)^n} \left\langle e^{i(x-y) \cdot \xi}, \, 
\phi(\varepsilon p(\xi)) \right\rangle_{\xi}    \, .  
\label{(4.4)}   
\end{equation}  
  
When $x \neq y$ are fixed,
  $e^{i(x-y) \cdot \xi}$ is distributionally small  
as a function of $\xi$. 
   This also holds distributionally in $(x,y) $.
   Thus  
the expansion of (\ref{(4.4)}) follows from the following lemma.  
  
\smallskip{\bf Lemma 4.1.} 
{\em Let $f \in {\cal K}'(\RR^n)$, 
  so that it satisfies the moment asymptotic  
expansion  
\begin{equation} f(\lambda x) \sim \sum_{k \in \NN^n} \frac{(-
1)^{|k|}\,\mu_k \,
  D^k  \delta(x)}  
{k! \; \lambda^{|k|+n}} \quad \mbox{as} \ \lambda \rightarrow   
\infty    ,  
\label{(4.5)}   
\end{equation}  
where $\mu_k=\langle f(x),x^k\rangle$ ,  $k \in \NN^n$, are the moments.  
Then  
if $p$ is an elliptic polynomial and $\phi \in {\cal K}$,
\begin{equation} \left\langle f(x),\,\phi(\varepsilon p(x))
  \right\rangle \sim   
\sum_{n=0}^{\infty} \frac{\langle f(x),\,p(x)^m \rangle\, \phi^{(m)}(0)   
\varepsilon^m}{m!}  \quad
\mbox{as} \ \varepsilon \rightarrow 0    .  
\label{(4.6)}   
\end{equation} 
 }
 \smallskip

{\em Proof:\/} The proof consists in showing that the Taylor expansion in   
$\varepsilon$,  
\begin{equation} \phi \left( \varepsilon p(x) \right) = \sum_{n=0}^{N}   
\frac{\phi^{(m)}(0) \,p(x)^m \varepsilon^m}{m!} + O(\varepsilon^{N+1})    ,  
\label{(4.7)}   
\end{equation}  
not only holds pointwise 
 but actually  holds in the topology of ${\cal K}(\RR^n)$.  
But the remainder in this Taylor approximation is
\begin{equation}  
R_N(x,\varepsilon)={\phi^{(N+1)}(\theta p(x))\,
 p(x)^{N+1}\varepsilon^{N+1}\over (N+1)!}\,,  
\label{(4.8)}  
\end{equation}  
for some $\theta\in (0,\varepsilon)$. 
 Since $\phi\in{\cal K}$, there exists $q\in\RR$ such  
that $\phi^{(j)}(x)=O(|x|^{q-j})$ as $x\to\infty$. 
 If $p$ has degree $m$ it  
follows that   
\begin{equation}  
|R_N(x,\varepsilon)|\leq {M\max{\{ 1,|x|^{mq}\} }\,\varepsilon^{N+1} 
 \over (N+1)!}   
\label{(4.9)}  
\end{equation}  
for some constant $M$, and the convergence of the Taylor expansion in the  
topology of the space ${\cal K}$ follows. $\;\;\bull$   
\smallskip

Thus, applying (\ref{(4.6)}) with $f(x)=e^{i(x-y) \cdot \xi}$
  for $x \neq y$ or  
 distributionally  
in $(x,y)$, we obtain  
\begin{displaymath}  
\Phi(\varepsilon)  \sim  \frac{1}{(2 \pi)^n} \sum_{k=0}^{\infty} \frac{  
\langle e^{i(x-y)\cdot \xi},\,p(\xi)^k \rangle \phi^{(k)}(0)  
 \varepsilon^k}{k!} \ ,  
\end{displaymath}  
or  
\begin{equation}  
\Phi (\varepsilon)  \sim  \sum_{k=0}^{\infty} \frac{H^k \delta(x-y)  
\, \phi^{(k)}(0) \varepsilon^k}  
{k!}\,.  
\label{(4.10)}  
\end{equation}  
Therefore,  
\begin{equation} e(x,y;\lambda \sigma) \sim \sum_{k=0}^{\infty}   
\frac{(-1)^k H^k \delta(x-y) \,\delta^{(k)}(\lambda)}{k! \,\sigma^{k+1}}  
 \quad  
\mbox{as} \ \sigma \rightarrow \infty    ,  
\label{(4.11)}   
\end{equation}  
in accordance with the general result.  
  
Observe also that if $H_1$ is any operator  
 corresponding to the same differential  
expression, 
 considered in some open set ${\cal M}$ with some boundary conditions,  
 then   
its spectral density $e_1(x,y;\lambda)$ satisfies  
\begin{equation} e_1(x,y;\lambda) = 
 \frac{1}{(2 \pi)^n} \left\langle e^{i(x-y)\cdot \xi},  \,
\delta(p(\xi)-\lambda)\right\rangle + o(\lambda^{\infty}) \quad 
(C)\quad
{\rm as}\ \lambda\to\infty.  
\label{(4.12)}   
\end{equation}  
        
 \smallskip
 {\bf Example.} Let ${\cal M}$ be a region in $\RR^n$ and let $H$ be any  
 self-adjoint  
extension of the negative Laplacian $-\Delta$ obtained by imposing suitable  
boundary conditions on ${\cal M}$. 
 Let $e_{{\cal M}}(x,y;\lambda)$ be the  
spectral density. Then  
\begin{equation}  
e_{{\cal M}}(x,y;\lambda)={1\over (2\pi)^n}\langle \delta (|\xi |^2-\lambda ),  
\,e^{i(x-y)\cdot\xi}\rangle + o(\lambda^{-\infty})\quad (C).  
\label{(4.13)}  
\end{equation}  
We now use the one-variable formula   
\begin{displaymath}  
\delta (f(x))={\delta (x-x_0)\over |f'(x_0)|}\, ,  
\end{displaymath}  
valid if $f$ has a single zero at $x_0\,$, and pass to polar coordinates  
$\xi=r\omega$, where $r=|\xi |$, 
$\omega=(\omega_1,\ldots ,\omega_n)$ satisfies  
$|\omega |=1$, and $d\xi = r^{n-1}\,dr\,d\sigma   
(\omega)$, to obtain  
\begin{eqnarray*}  
{1\over (2\pi)^n}\langle \delta (|\xi |^2-\lambda),\,
 e^{i(x-y)\cdot\xi}\rangle  
& = & {1\over (2\pi)^n}\int_{|\omega |=1}
 \int_0^\infty \langle \delta (r^2-  
\lambda),\,e^{ir(x-y)\cdot \omega}\rangle r^{n-1}\,
  dr\,d\sigma (\omega) \\  
& = & {\lambda^{n/2 -1}\over 2(2\pi)^n}\int_{|\omega |=1}e^{i\lambda^{1/2}  
(x-y)\cdot\omega}\,d\sigma (\omega) \\  
& = & {\lambda^{n/2 -1}\over 2(2\pi)^n}\int_{|\omega |=1}e^{i\lambda^{1/2}  
\omega_1|x-y|}\,d\sigma (\omega ) \\  
& = & {\lambda^{n/2 -1}\over 2(2\pi)^n}{2\pi^{(n-1)/2}\over  
 \Gamma ({n-1\over 2})}  
\int_{-1}^1e^{i\lambda^{1/2}u|x-y|}(1-u^2)^{{n-3\over 2}}\,du \\  
& = & {\lambda^{n/4 -1/2}J_{n/2 -1}(\lambda^{1/2}|x-y|)\over  
 2^{n/2 +1}\pi^{n/2}  
|x-y|^{n/2 -1}}\, ,  
\end{eqnarray*}  
where $J_p(x)$ is the Bessel function of order $p$. Therefore  
\begin{equation}  
e_{{\cal M}}(x,y;\lambda)= {\lambda^{n/4 -1/2}J_{n/2 -1}(\lambda^{1/2}  
|x-y|)\over 2^{n/2 +1}\pi^{n/2} |x-y|^{n/2 -1}}+  
 o(\lambda^{-\infty}) \quad (C)\quad {\rm as} 
\ \lambda\to\infty,  
\label{(4.14)}  
\end{equation}  
uniformly over compacts of ${\cal M}\times{\cal M}$. $\;\;\bull$  
\smallskip   

 Our second model is an ordinary differential operator $H$ with  
variable coefficients, as treated in Refs.~\onlinecite{F1,F2,F3}.  
 There are two major simplifications in this one-dimensional case. 
 First, the Weyl--Titchmarsh--Kodaira theory\cite{T,Ko} 
 expresses the spectral density as 
 \begin{equation}e(x,y;\lambda)\, d\lambda = 
 \sum_{j,k=0}^1 \psi_{\lambda j}(x) \, d\mu^{jk}(\lambda) \,  
 \overline{\psi_{\lambda k}(y)}\,, \label{(4.15)} 
 \end{equation} 
 where $d\mu^{jk}$ are certain Stieltjes measures supported on the  
spectrum of $H$, and $\psi_{\lambda j}$ are the classical solutions  
of $H\psi - \lambda \psi$ with the basic data 
  \begin{eqnarray} \psi_{\lambda0} (x_0) &=&1,
  \quad \psi'_{\lambda0} (x_0) = 0, \nonumber\\
 \psi_{\lambda1} (x_0) &=&0, \quad \psi'_{\lambda1} (x_0) = 1,
 \label{(4.16)} \end{eqnarray} 
 at some $x_0 \in {\cal M}$. 
 Thus  
 \begin{eqnarray} \mu^{00}(\lambda) &=& E(x_0,x_0;\lambda), \quad 
\mu^{01}(\lambda) = {\partial E\over \partial y}(x_0,x_0;\lambda), 
\nonumber\\ 
\mu^{10}(\lambda) &=& {\partial E\over \partial x}(x_0,x_0;\lambda), 
\quad 
  \mu^{11}(\lambda) = {\partial^2 E\over \partial x\,\partial y} 
 (x_0,x_0;\lambda). 
  \label{(4.17)}\end{eqnarray} 
 Second, the eigenfunctions $\psi_{\lambda j}$ can be approximated  
for large $\lambda$ quite explicitly by the phase-integral (WKB)  
method. 
 (Thirdly, but less essentially, there is no loss of generality in  
considering 
 \begin{equation}H = - \, {d^2 \over dx^2} + V(x),
  \label {(4.18)}\end{equation} 
 since the general second-order operator can be reduced to this  
form by change of variables.) 
 
 In Ref.~\onlinecite{F1}  the phase-integral representation of the  
eigenfunctions was used to obtain in a direct and elementary way  
the expansion 
 \begin{equation}d\mu^{jk}(\lambda) \sim 
 {1\over \pi} \sum_{n=0}^\infty  \rho_n^{jk}(x_0)  
\, \omega^{2 \delta_{j1}\delta_{k1} - 2n} \, d \omega, 
 \label{(4.19)}\end{equation} 
 where $\lambda = \omega^2$ and 
 \begin{eqnarray}\rho_0^{00}&=&1, \quad \rho_1^{00}= {1\over2} V, \quad    
 \rho_2^{00}={1\over 8} (-V'' + 3V^2), \quad\ldots, 
 \nonumber\\ 
\rho_0^{11}&=&1, \quad \rho_1^{11}= -\,{1\over2} V, \quad    
 \rho_2^{11}={1\over 8} (V'' - 3V^2), \quad\ldots, \label {(4.20)}
 \end{eqnarray} 
 \begin{displaymath}\rho_n^{10} =\rho_n^{01} = {1\over2}\, {d\over dx_0}  
(\rho_n^{00}).\end{displaymath} 
 
 Formula (\ref{(4.19)}) is a rigorous asymptotic expansion when 
 ${\cal M} =  
\RR$ and $V$ is a $C^\infty$ function of compact support. 
 The relevance of (\ref{(4.19)}) in more general cases, where it is  
certainly not a literal pointwise asymptotic expansion, 
 was discussed at length in Ref.~\onlinecite{F1}; 
 the results of the present paper simplify and sharpen that  
discussion by showing that 
 {\em the error in (\ref{(4.19)}) is $O(\lambda^{-\infty})$ in the $(C)$  
sense for any operator locally equivalent to one for which
 (\ref{(4.19)})  
holds pointwise.}

\section{Pointwise and average expansions}  
\label{sec5}

Let $f \in {\cal K}'(\RR)$.   
 Since the elements of ${\cal K}'(\RR)$ are precisely the distributionally  
small generalized functions,  
 it follows that $f$ satisfies the moment asymptotic  
expansion; that is,  
\begin{equation} f(\lambda x) \sim 
 \sum_{j=0}^{\infty} \frac{(-1)^j \, \mu_j \, \delta^{(j)}  
(x)}{j! \, \lambda^{j+1}}   \quad \mbox{as} \ \lambda \rightarrow  
\infty    ,  
\label{(5.1)}   
\end{equation}  
where  
\begin{equation} \mu_j = \langle f(x), x^j \rangle , \quad j \in \NN ,  
\label{(5.2)}   
\end{equation}  
are the moments.  
  
The moment asymptotic expansion allows us to obtain the 
 small-$t$ behavior of  
functions $G(t)$ that can be written as  
\begin{equation} G(t) = \langle f(x), g(tx) \rangle  ,   
\label{(5.3)}   
\end{equation}  
as long as $g \in {\cal K}$.   Indeed, (\ref{(5.1)}) gives  
\begin{equation} G(t) = \sum_{j=0}^{\infty} \frac{\mu_j
  \, g^{(j)}(0) \, t^j}{j!}   
 \quad   
\mbox{as} \ t \rightarrow 0    .  
\label{(5.4)}   
\end{equation}  
  
Naturally, this would be valuable if $f(\lambda)=e(x,y;\lambda)$ 
is the spectral  
density of the elliptic differential operator $H$ and $G(t,x,y)=\langle   
e(x,y;\lambda),\, g(\lambda t) \rangle$ is an associated Green kernel.  
  
However, we emphasize that the derivation of  (\ref{(5.4)})
  holds only 
 when $g \in {\cal K}$. 
   What if   
$g \notin {\cal K}\,$? 
   A particularly interesting example is the kernel $U(t,x,y)=  
\langle e(x,y;\lambda),\, e^{-i \lambda t} \rangle$ 
that solves the Schr\"odinger  
equation  
 \begin{mathletters}\label{(5.5)}
\begin{equation}i \frac{\partial U}{\partial t} =  HU  , \quad t>0  
\label{(5.5a)}   
\end{equation}  
with initial condition  
\begin{equation} U(0^{+},x,y) = \delta(x-y)    .  
\label{(5.5b)}   
\end{equation}  \end{mathletters}
In this case $g(x)=e^{-ix}$ is smooth,
  but because of its behavior at infinity,
it does not belong to ${\cal K}$.  
   We pointed out in the introduction, however, that
(\ref{(5.4)}) is still valid in some ``averaged'' sense.  
  
Indeed, we  
 shall now show that formula (\ref{(5.3)})
  permits one to define $G(t)$  
 as a distribution  
when instead of asking $g \in {\cal K}$ we assume
  $g$ to be a tempered distribution  
of the space ${\cal S}'$
  which has a distributional expansion at the origin. 
   We then  
show that (\ref{(5.4)}) holds in an averaged or distributional sense.  
 The fact that  
the space of smooth functions ${\cal K}$  
 is replaced by the space of tempered distributions is  
not casual: 
the distributions of ${\cal S}'$ are exactly those that have the behavior  
at $\infty$ of the elements of ${\cal K}$ 
 in the Ces\`aro or distributional sense.    
Indeed, we have  
 
\smallskip{\bf Lemma 5.1} 
{\em  Let $g\in {\cal S}'(\RR)$.   Then there exists $\alpha \in \RR$  
 such that  
\begin{equation} g^{(n)}(\lambda x) 
 = O(\lambda^{\alpha-n})   \quad \mbox{as} \
\lambda \rightarrow \infty    ,  
\label{(5.6)}   
\end{equation}  
distributionally.
}
\smallskip  
  
{\em Proof:\/} See Ref.~\onlinecite{E96}, where it is shown that 
 (\ref{(5.6)}) is actually   
a characterization of the tempered distributions. $\;\; \bull$  
\smallskip

Let $g \in {\cal S}'(\RR)$ and let $\alpha$ be as in (\ref{(5.6)}).   
 If $\phi \in {\cal S}(\RR)$ then the  
function $\Phi$ defined by  
\begin{equation} \Phi(x) = \langle g(tx),\phi(t) \rangle  
\label{(5.7)}   
\end{equation}  
is smooth in the open set $(-\infty,0) \cup (0,\infty)$ and,  
 because of (\ref{(5.6)}), satisfies  
\begin{equation} \Phi^{(n)}(x) = 
 O(|x|^{\alpha-n}) \quad\mbox{as} 
  \ |x| \rightarrow  \infty    .  
\label{(5.8)}   
\end{equation}  
It follows that we can define $G(t)=\langle f(x),g(tx) \rangle$ 
 as an element  of ${\cal S}'(\RR)$ by  
\begin{equation} \langle G(t),\phi(t) \rangle =
  \langle f(x),\Phi(x) \rangle    ,  
\label{(5.9)}   
\end{equation}  
whenever $f \in {\cal K}'$ and $0 \notin \mathop{\rm supp} f$.  
  
When $0 \in \mathop{\rm supp} f$ then (\ref{(5.9)})
  cannot be used unless $\Phi$ is smooth at the  origin.  
  And in order to have $\Phi$ smooth we need to ask the existence of   
the {\em distributional\/} values $g^{(n)}(0)$,  $n=0,1,2, \ldots$.  
   
Recall that following {\L}ojasiewicz,\cite{L}  one says that a 
distribution  $h \in {\cal D}'$  
has the value $\gamma$ at the point $x=x_0\,$, written as  
\begin{equation} h(x_0) = \gamma  \quad \mbox{in} \ {\cal D}',  
\label{(5.10)}   
\end{equation}  
if  
\begin{equation} \lim_{\varepsilon \rightarrow 0} h(x_0+\varepsilon x) 
 = \gamma   
\label{(5.11)}   
\end{equation}  
distributionally; that is, if for each $\phi \in {\cal D}$  
\begin{equation} \lim_{\varepsilon \rightarrow 0} \,
 \langle h(x_0+\varepsilon x) , \phi(x) \rangle =  
\gamma \int_{-\infty}^{\infty} \phi(x)\, dx    .  
\label{(5.12)}   
\end{equation}  
It can be shown that $h(x_0)=\gamma$ in ${\cal D}'$ 
  if and only if there exists a primitive  
$h_n$ of some order~$n$, $ h_n^{(n)}=h$,
  which is continuous in a neighborhood  
of $x=x_0$ and satisfies  
\begin{equation} h_n(x) = \frac{\gamma(x-x_0)^n}{n!} + o(|x-x_0|^n) ,  
 \quad \mbox{as} \
x \rightarrow x_0   \, .  
\label{(5.13)}   
\end{equation}  
  
In our present case, we need to ask the existence of the 
 distributional values  
$g^{(n)}(0)=a_n$ for $n \in \NN$.   
  We can then say that $g(x)$ has the small-$x$  
``averaged'' or distributional expansion  
\begin{equation} g(x) \sim \sum_{n=0}^{\infty} \frac{a_n \,x^n}{n!} 
\, ,   \quad \mbox{as} \  
x \rightarrow 0 , \quad \mbox{in} \ {\cal D}'    ,  
\label{(5.14)}   
\end{equation}  
in the sense that the parametric expansion  
\begin{equation} g(\varepsilon x) \sim \sum_{n=0}^{\infty} 
 \frac{a_n \, \varepsilon ^n \ x^n}{n!}   
\, , \quad \mbox{as} \ \varepsilon  \rightarrow 0     ,  
\label{(5.15)}   
\end{equation}  
holds, or, equivalently, that   
\begin{equation} \langle g(\varepsilon x),\phi(x) \rangle
  \sim \sum_{n=0}^{\infty} \frac{a_n}{n!}   
\left( \int_{-\infty}^{\infty} x^n \phi(x)\,dx \right) \varepsilon ^n,   
\label{(5.16)}   
\end{equation}  
for each $\phi \in {\cal D}$.  
  
\smallskip{\bf Lemma 5.2.}
{\em  Let $g \in {\cal S}'$ be such that the distributional values  
 $g^{(n)}(0)=a_n \, ,  \mbox{in}\  {\cal D}'$, exist for $n \in \NN.$  
Let $\phi \in {\cal S}$ and  
put $\Phi(x)=\langle g(tx),\phi(t)\rangle. $   
 Then $\Phi \in {\cal K}$. 
} 
\smallskip  

{\em Proof:\/} Indeed, $\Phi$ is smooth for $x\not= 0$, but since the 
distributional values $g^{(n)}(0)=a_n$ exist, it follows that 
$\Phi (x)\sim\sum_{n=0}^\infty b_nx^n$ as $x\to 0$, where
$b_n=(a_n/n!)\int_{-\infty}^\infty x^n\phi (x)\,dx$. 
 Thus $\Phi$ is
also smooth at $x=0$. 
 Finally, let $\alpha$ be as in (\ref{(5.6)});
  then
$\Phi^{(n)}(x)=O(|x|^{\alpha -n})$ as $|x|\to\infty$. Hence
$\Phi\in{\cal K}$. $\;\; \bull$  
\smallskip  
  
 Using this lemma we can give the following  
  
 \smallskip
{\bf Definition.}
{\em Let $f \in {\cal K}'$.   
  Let $g \in {\cal S}'$ have distributional values  
$g^{(n)}(0)$,  $n \in \NN$. 
  Then we can define the tempered distribution  
\begin{equation} G(t) = \langle f(x),g(tx) \rangle  
\label{(5.17)}   
\end{equation}  
by  
\begin{equation} 
 \langle G(t),\phi(t) \rangle = \langle f(x),\Phi(x) \rangle  ,  
\label{(5.18)}   
\end{equation}  
where  
\begin{equation} \Phi(x) = \langle g(tx),\phi(t) \rangle    ,  
\label{(5.19)}   
\end{equation}  
if $\phi \in {\cal S}.$
}  
\smallskip

In general the distribution $G(t)$ is not smooth near the origin, 
  but its distributional  
behavior can be obtained from the moment asymptotic expansion.  

\smallskip
{\bf Theorem 5.1.}  
{\em Let $f \in {\cal K}'$ with moments $\mu_n=\langle f(x),x^n \rangle$.     
Let $g \in {\cal S}'$ have distributional values $g^{(n)}(0)$ 
 for $n \in \NN$. 
   Then  
the tempered distribution $G(t)=\langle f(x),g(tx) \rangle$ 
has distributional   
values $G^{(n)}(0)$, $ n \in \NN$,
  which are given by $G^{(n)}(0)=\mu_n g^{(n)}  (0)$, 
 and $G$ has the distributional expansion  
\begin{equation} G(t) \sim \sum_{n=0}^{\infty} 
 \frac{\mu_n g^{(n)}(0) t^n}{n!} \, ,  
 \quad \mbox{in} \
{\cal D}' , \quad \mbox{as} \ t \rightarrow 0    .  
\label{(5.20)}   
\end{equation}
}
\smallskip
  
 {\em Proof:\/}   
 Let $\phi \in {\cal S}$ and let $\Phi(x)=\langle g(tx),\phi(t) \rangle$.   
Then  
\begin{equation} \langle G(\varepsilon t),\phi(t) \rangle = 
 \langle f(x),\Phi(\varepsilon x)\rangle ,  
\label{(5.21)}   
\end{equation}  
and since $\Phi^{(n)}(0)=g^{(n)}(0) \int_{-\infty}^{\infty}t^n \phi(t)\,dt$,
  the  moment asymptotic expansion yields  
\begin{equation} \langle G(\varepsilon t),\phi(t) \rangle 
  \sim \sum_{n=0}^{\infty}   
\frac{\mu_n \, g^{(n)}(0)}{n!} \left( \int_{-\infty}^{\infty} t^n 
\phi(t)\,dt   
\right) \varepsilon ^n   \quad \mbox{as} \ \varepsilon \rightarrow 0 ,   
\label{(5.22)}   
\end{equation}  
and (\ref{(5.20)}) follows. $\;\; \bull $  
\smallskip

Before we continue, it is worthwhile to give some examples.  

 \smallskip  
{\bf Example.}
  Let $g \in {\cal S}'$ be such that the distributional values $g^{(n)}  
(0)$ exist for $n \in \NN$.   
Since the Fourier transform $\hat{g}(\lambda)$ can  
be written as 
 $\hat{g}(\lambda)=\lambda^{-1} \langle e^{ix},g(\lambda^{-1} x)  
 \rangle$, 
 and since all the moments $\mu_n= \langle e^{ix},x^n\rangle$ vanish,
it follows that 
 $\hat{g}(\varepsilon^{-1})=O(\varepsilon^{\infty})$ distributionally  
as $\varepsilon \rightarrow 0$ 
 and thus $\hat{g}(\lambda)=O(|\lambda|^{-\infty}) \ \
(C) \ \ \mbox{as} \ |\lambda| \rightarrow \infty. $   
 Therefore $\hat{g}  
\in {\cal K}'$.  
  
Conversely, if $f \in {\cal K}'$, then its Fourier transform 
 $\hat{f}(t)$ is equal to
$F(t)=\langle f(x),e^{itx}\rangle$ for $t \neq 0$. 
   Thus $\hat{f}(t)=F(t)+  
\sum_{j=0}^{n} a_j \, \delta^{(j)}(t)$  
 for some constants $a_0\,, \ldots, a_n\,$.     
But the distributional values $F^{(n)}(0)$ exist for $n \in \NN$ 
 and are given   
by $F^{(n)}(0)=i^n  \langle f(x),x^n \rangle$,  
 and hence $\hat{F} \in {\cal K}'$, and  
it follows that $a_0=\cdots=a_n=0$. 
   In summary, $\hat{f}^{(n)}(0)$ exists in  
${\cal D}'$ for each $n \in \NN$.  
  
Therefore, a distribution $g \in {\cal S}'$ is smooth at the origin  
 in the distributional  
sense 
 (that is, the distributional values $g^{(n)}(0)$ exist for $n \in \NN$) 
 if and only if its Fourier transform $\hat{g}$ is distributionally small 
 (i.e.,  $\hat{g} \in {\cal K}'$). $\;\; \bull$  
 \smallskip
  
{\bf Example.} Let $\xi \in \CC$ with $|\xi|=1$, $\xi \neq 1$.    
 Then the distribution  
$f(x)=\sum_{n=-\infty}^{\infty} \xi^n \delta(x-n)$ 
 belongs to ${\cal K}'$.
    All   the moments vanish:  
 $\mu_k=\sum_{n=-\infty}^{\infty} \xi^n n^k=0 \ \ (C) \ \
\mbox{for} \ k=0,1,2,\ldots$.    
 It follows that if $g \in {\cal S}'$ is distributionally  
smooth at the origin, then  
\begin{equation} \sum_{n=-\infty}^{\infty} \xi^n g(nx) = o(x^{\infty}) 
  \quad \mbox{in} \
{\cal D}' \quad \mbox{as} \ x \rightarrow 0    .  
\label{(5.23)}   
\end{equation}  
  
When $\xi=1$,  $\sum_{n=-\infty}^{\infty} \delta(x-n)$ 
does not belong to ${\cal K}'$  
but $\sum_{n=-\infty}^{\infty} \delta(x-n)-1$ does. 
  Thus, if $g \in {\cal S}'$ is  
distributionally smooth at the origin and
  $\int_{-\infty}^{\infty} g(u)\,du$ 
is  defined, then  
\begin{equation} \sum_{n=-\infty}^{\infty} g(nx) =
  \left( \int_{-\infty}^{\infty} 
g(u)\,du \right)  
x^{-1} + o(x^{\infty}) \quad \mbox{in} \ {\cal D}'   
\quad \mbox{as} \   x \rightarrow 0    .  
\label{(5.24)}   
\end{equation}  
Actually, many number-theoretical expansions  considered in  
Ref.~\onlinecite{E92} and Chapter~5 of Ref.~\onlinecite{EK94}
 will hold in the averaged or distributional sense when   
applied to distributions. $\;\; \bull$  
\smallskip

Many times, 
  $\mathop{\rm supp} f \subseteq [0,\infty)$ 
 and one is interested in $G(t)=\langle  
f(x),g(tx) \rangle$ for $t>0$ only. 
  In those cases the values of $g(x)$ for  
$x<0$ are irrelevant and one may assume that
  $\mathop{\rm supp} g \subseteq [0,\infty)$.  
Since we need to consider $\Phi(x)=\langle g(tx),\phi(t) \rangle$ for $x>0$ 
only,  we do not require the existence of the distributional values 
$g^{(n)}(0)$;   instead, we assume the existence of the one-sided 
distributional values $g^{(n)}(0^{+})  =a_n$ for $n \in \NN$.
    This is 
equivalent to asking $g(\varepsilon x)$ to have the  asymptotic development 
 \begin{equation} g(\varepsilon x) \sim
  \sum_{n=0}^{\infty} \frac{a_n \,\varepsilon^n \,x_{+}^n}{n!}    \quad 
\mbox{as} \ \varepsilon \rightarrow 0^{+} \, ;   \label{(5.25)} 
   \end{equation}  
that is,
  \begin{equation} 
\langle g(\varepsilon x),\phi(x) \rangle  \sim \sum_{n=0}^{\infty}
    \frac{a_n}{n!} 
\left( \int_{0}^{\infty} x^n \phi(x)\,dx \right) \varepsilon ^n 
\quad
\mbox{as} \ \varepsilon \rightarrow 0^{+}       \label{(5.26)} 
   \end{equation} 
 for $\phi \in {\cal S}$.
    We shall use the notation  
\begin{equation} g(x) \sim \sum_{n=0}^{\infty} \frac{a_n 
x^n}{n!} \quad \mbox{in}   \ {\cal D}' 
  \quad \mbox{as} \ x 
\rightarrow 0^{+}
    \label{(5.27)}   \end{equation} 
 in such a case.  
  
\smallskip{\bf Lemma 5.3.} 
{\em  Let $g \in {\cal S}'$ with 
 $\mathop{\rm supp} g \subseteq [0,\infty)$ and let
  $g^{(n)}(0^{+})=a_n$ exist in ${\cal D}'$ for $n \in \NN$.
    Let $\phi \in {\cal S}$ 
and   put $\Phi(x)=\langle g(tx),\phi(t) \rangle$ for $ x>0$. 
   Then $\Phi$  
admits extensions $\tilde{\Phi}$ to $\RR$ with 
 $\tilde{\Phi} \in {\cal K}(\RR).$
}
\smallskip

{\em Proof:\/} It suffices to show that $\Phi$ is smooth up to the
origin from the right and that it satisfies estimates of the form
$\Phi^{(j)}(x)=O(|x|^{\alpha -n})$ as $|x|\to\infty$. 
 But the first 
statement follows because $g^{(n)}(0^+)$ exists for all $n\in\NN$,
while the latter is true because of (\ref{(5.6)}). $\;\;\bull$ 
 \smallskip  

 From this lemma it follows that when $f \in {\cal K}'$, 
 $\mathop{\rm supp} f \subseteq   
[0,\infty)$, $\mathop{\rm supp} g \subseteq [0,\infty)$,
  and the distributional values 
$g^{(n)}  (0^{+})$ exist for $n \in \NN$, 
 then $G(t)=\langle f(x),g(tx) \rangle$ 
 can be  defined as a tempered distribution with support contained in 
$[0,\infty)$ by  
 \begin{equation} \langle G(t),\phi(t) \rangle = \langle 
f(x),\tilde{\Phi}(x) \rangle    , 
  \label{(5.28)}   \end{equation} 
  where $\tilde{\Phi}$ is 
any extension of $\Phi(x)=\langle g(tx),\phi(t) \rangle $, $ x>0$, 
 such that $\tilde{\Phi} \in {\cal K}$.  
  
\smallskip{\bf Theorem 5.2.} 
{\em  Let $f \in {\cal K}'$ with $\mathop{\rm supp} f \subseteq [0,\infty)$ 
and moments  $\mu_n=\langle f(x),x^n \rangle$.
    Let $g \in {\cal S}'$ with $\mathop{\rm supp} g 
\subseteq [0,\infty)$   have distributional one-sided values 
 $g^{(n)}(0^{+})$ for $n \in \NN$. 
   Then the tempered  distribution $G(t)=\langle f(t),g(tx) 
\rangle$ defined by (\ref{(5.28)}) has distributional 
  one-sided values 
$G^{(n)}(0^{+})$, $n \in \NN$, which are given by $G^{(n)}(0^{+})=  \mu_n 
\,g^{(n)}(0^{+})$, and $G$ has the distributional expansion
   \begin{equation} G(t) \sim 
\sum_{n=0}^{\infty} \frac{\mu_n \,g^{(n)}(0^{+}) \,t^n}{n!} \quad
\mbox{in} \ {\cal D}'     \quad \mbox{as} \ t \rightarrow 0^{+}    .  
\label{(5.29)}   \end{equation}
}
\smallskip  
   
{\em Proof:\/} Quite similar to the proof of Theorem 5.1 . $\;\; \bull$  
\smallskip

\section{Expansion of Green kernels I{}: Local expansions}  
 \label{sec6}

In this section we shall consider the small-$t$ behavior of Green kernels of  
the type $G(t;x,y)=\langle e(x,y;\lambda),g(\lambda t)\rangle$ 
for some $g \in   {\cal S}'. $  
 Here $e(x,y;\lambda)$ is the spectral density kernel corresponding to  
a positive elliptic operator $H$ that acts on the smooth 
 manifold~${\cal M}$.  
  
Our results can be formulated in a general framework.   
 So, let $H$ be a positive  
self-adjoint operator on the domain ${\cal X}$ of the 
 Hilbert space~${\cal H}$.   
  Let ${\cal X}_{\infty}$  
be the common domain of $H^n$,  $n \in \NN$, 
  and let $e_{\lambda}$ be the associated spectral  
density.  
 Let $g \in {\cal S}'$ with $\mathop{\rm supp} g  \subseteq [0,\infty)$ 
 such that the  
one-sided distributional values $a_n=g^{(n)}(0^{+})$ exist for 
 $n \in \NN$.     
Then we can define  
\begin{equation} G(t) = g(tH)    , \quad t>0    ,  
\label{(6.1)}   
\end{equation}  
that is,  
\begin{equation} G(t) = \langle e_{\lambda},g(t \lambda)\rangle  ,
  \quad t>0 \; .  
\label{(6.2)}   
\end{equation}  
Thus $G$ can be considered an operator-valued distribution  
 in the space ${\cal S}'(\RR,  L({\cal X}_{\infty},{\cal H}))$.
    The behavior of $G(t)$ as $t \rightarrow 0^{+}$   
can be obtained from the moment asymptotic expansion (\ref{(3.9)})
 for $e_{\lambda}\,$.
  The   
expansion of $G(t)$ as $t \rightarrow 0^{+}$ will be a distributional or   
``averaged'' expansion, in general,  
 but when $g$ has the behavior of the elements  
of ${\cal K}$ at $\infty$ it becomes a pointwise expansion.   
 In particular, if $g$ is  
smooth in $[0,\infty)$,  
 the expansion is pointwise or not depending on the behavior  
of $g$ at infinity.

\smallskip{\bf Theorem 6.1.}  
{\em Let $H$ be a positive self-adjoint operator on the domain  
${\cal X}$ of the Hilbert space~${\cal H}$.
    Let ${\cal X}_{\infty}$ be the intersection of the  
domains of $H^n$ for $n \in \NN$.   
  Let $g \in {\cal S}'$ with $\mathop{\rm supp} g \subseteq [0,\infty)$  
be such that the distributional one-sided values  
\begin{equation} g^{(n)}(0^{+}) = a_n   
  \quad \mbox{in} \ {\cal D}'    
\label{(6.3)}   
\end{equation}  
exist for $n \in \NN$. 
   Let $G(t)=g(tH)$, an element of ${\cal S}'(\RR,L({\cal X}_{\infty},  
{\cal H}))$ with support contained in $[0,\infty)$.    
 Then $G(t)$ admits the distributional  
expansion in $L({\cal X}_{\infty},{\cal H})$,  
\begin{equation} G(t) \sim \sum_{n=0}^{\infty} \frac{a_n H^n t^n}{n!}   
\, , \quad \mbox{as} \
t \rightarrow 0^{+}  , \quad \mbox{in} \ {\cal D}'    ,  
\label{(6.4)}   
\end{equation}  
so that the distributional one-sided values $G^{(n)}(0^{+})$  
 exist and are given  
by  
\begin{equation} G^{(n)}(0^{+}) = a_n H^n   
 \quad \mbox{in} \ {\cal D}'    .  
\label{(6.5)}   
\end{equation}  
When $g$ admits an extension that belongs to ${\cal K}$,  
 (\ref{(6.4)}) is an ordinary pointwise  
expansion while the $G^{(n)}(0^{+})$ exist as ordinary one-sided values.
}
\smallskip  
  
{\em Proof:\/} Follows immediately from Theorem 5.2.
  $\;\; \bull$  
\smallskip

When $H$ is a positive elliptic differential operator  
 acting on the manifold ${\cal M}$,   
then  Theorem~6.1 gives  the small-$t$
  expansion of Green kernels.  Let  
$e(x,y;\lambda)$ be the spectral density kernel and let  
\begin{equation} G(t,x,y) = \langle e(x,y;\lambda),g(t \lambda)\rangle    ,
  \quad t>0,   
\label{(6.6)}   
\end{equation}  
be the Green function kernel corresponding to the operator $G(t)=g(tH)$.   
Then  
$G$ belongs to ${\cal S}'(\RR) \hat{\otimes}
  {\cal D}'({\cal M} \times {\cal M} )$, 
  has spectrum in $[0,\infty)$,  
and as $t \rightarrow 0^{+}$ admits the distributional expansion  
\begin{equation} G(t,x,y) \sim \sum_{n=0}^{\infty}
  \frac{a_n \,H^n \delta(x-y) \,t^n}{n!} \, , 
  \quad   
\mbox{as} \ t \rightarrow 0^{+}    , \quad \mbox{in} \ {\cal D}' ;
\label{(6.7)}   
\end{equation}  
that is,  
\begin{equation} G(\varepsilon t,x,y) \sim \sum_{n=0}^{\infty}
  \frac{a_n \,H^n \delta(x-y) \,\varepsilon^n  \,
t_{+}^n}{n!}\,   \quad \mbox{as} \ \varepsilon  \rightarrow 0^{+} , 
\label{(6.8)}   
\end{equation}  
in ${\cal D}'({\cal M}\times{\cal M})$.  
  Also, the distributional one-sided values   
${\displaystyle
 \frac{\partial^n}{\partial t^n}G(0^{+},x,y)}$ exist for $n \in \NN$   
and are given by  
\begin{equation} \frac{\partial^n}{\partial t^n}G(0^{+},x,y) = 
 a_n\, H^n \delta(x-y) 
 \quad  
\mbox{in} \ {\cal D}'    .  
\label{(6.9)}   
\end{equation}  
If $g$ admits extension to ${\cal K}$,  
 then (\ref{(6.7)}) and (\ref{(6.9)}) are valid in the ordinary  
pointwise sense with respect to $t$ (and distributionally in $(x,y)$).  
  
Pointwise expansions in $(x,y)$ follow when $x \neq y$.
    Indeed, if $U$ and  
$V$ are open subsets of ${\cal M}$ with  
 $U \cap V=\emptyset$, then $G$ belongs to ${\cal S}'  
(\RR) \hat{\otimes} {\cal E}(U\times V)$ 
 and as $t \rightarrow 0^{+}$ we have the   
distributional expansion  
\begin{equation} G(t,x,y) = o(t^{\infty})    ,
  \quad \mbox{in} \ {\cal D}'  , \quad
\mbox{as} \ t \rightarrow 0^{+},  
\label{(6.10)}   
\end{equation}  
in ${\cal E}(U\times V)$, and in particular pointwise on $x \in U$ 
and $y \in V$.
      The expansion becomes pointwise in $t$ when $g$ 
admits an extension to~${\cal K}$.  
  
These expansions depend only on the local behavior  
 of the differential operator.  
Let $H_1$ and $H_2$ be two differential operators that coincide on the open  
subset $U$ of~${\cal M}$.  
  Let $e_1(x,y,\lambda)$, $e_2(x,y,\lambda)$ be the corresponding  
spectral densities and $G_1(t,x,y)$ and $G_2(t,x,y)$ the 
 corresponding kernels  
for the operators $g(tH_1)$ and $g(tH_2)$, respectively. 
  Then  
\begin{equation} G_1(t,x,y) = G_2(t,x,y) + o(t^{\infty})      ,  
   \quad \mbox{in} \ {\cal D}'    , \quad
\mbox{as} \ t \rightarrow 0^{+}    ,  
\label{(6.11)}   
\end{equation}  
in ${\cal E}(U\times U)$; and when $g$ admits an extension  
that belongs to ${\cal K}$ this also holds pointwise in $t$.  
  
Let us consider some illustrations.

     \smallskip
{\bf Example.} 
  Let $K(t,x,y)=\langle e(x,y;\lambda),\,e^{-\lambda t}\rangle$   
be the heat kernel, corresponding to the operator $K(t)=e^{-tH}$, so that  
\begin{equation} \frac{\partial K}{\partial t} = -HK  , \quad t>0    ,  
\label{(6.12)}   
\end{equation}  
and  
\begin{equation} K(0^{+},x,y) = \delta(x-y)    .  
\label{(6.13)}   
\end{equation}  
In this case $g(t)=\chi (t)e^{-t}$ admits extensions in ${\cal K}$.  
  Thus the expansions  
\begin{equation} K(t,x,y) \sim \sum_{n=0}^{\infty} 
 \frac{(-1)^n\, H^n \delta(x-y)\,t^n}{n!} 
\quad\mbox{as} \ t \rightarrow 0^{+}  
\label{(6.14)}   
\end{equation}  
in the space ${\cal D}'({\cal M}\times{\cal M})$, and   
\begin{equation} K(t,x,y) = o(t^{\infty}) \quad \mbox{as} \ t  
\rightarrow 0^{+}, \quad \mbox{with} \ x \neq y, 
\label{(6.15)}   
\end{equation}  
  hold {\em pointwise\/} in $t$. $\;\;\bull$  
\smallskip

{\bf Example.}
  Let $U(t,x,y)=\langle e(x,y;\lambda),\,e^{-i \lambda t} \rangle$  
be the Schr\"odinger kernel, corresponding to $U(t)=e^{-itH}$, so that  
\begin{equation}i \frac{\partial U}{\partial t} = HU      , \quad t>0   
\label{(6.16)}   
\end{equation}  
and  
\begin{equation} U(0^{+},x,y) = \delta(x-y)    .   
\label{(6.17)}   
\end{equation}  
Here the function $e^{-it}$ belongs to ${\cal S}'$ but not to~${\cal K}$. 
   Therefore, the   expansions  
\begin{equation} U(t,x,y) \sim \sum_{n=0}^{\infty} 
 \frac{(-i)^n\, H^n \delta(x-y)\,t^n}{n!} 
  \, , \quad  
\mbox{as} \ t \rightarrow 0^{+}      , \quad \mbox{in} \ {\cal D}',   
\label{(6.18)}   
\end{equation}  
and  
\begin{equation} U(t,x,y) = o(t^{\infty})     ,
  \quad \mbox{in} \ {\cal D}' , 
 \quad   
\mbox{as} \ t \rightarrow 0^{+} \quad \mbox{with} \ x \neq y, 
\label{(6.19)}   
\end{equation}  
 are distributional or ``averaged'' in $t$.  
  
Consider, for instance, $U(t,x,y)=(4 \pi t)^{-n/2}e^{-i|x-y|^2/4t}$,   
 corresponding to the case when $H$   
is the negative Laplacian, $- \Delta$, acting on ${\cal M} =\RR^n$.    
 If $x \neq y$ are fixed,   
$U(t,x,y)$ oscillates as $t \rightarrow 0^{+}$, 
  but (\ref{(6.19)}) holds in the distributional  
sense. 
  (This also follows from the results of Ref.~\onlinecite{E96}, because   
$\displaystyle\lim_{s \rightarrow \infty} e^{ias}=0 \ \ (C) 
 \ \ \mbox{if} \ a \neq 0$.)  
  
Notice that it also follows that if $U_1$ is the Schr\"odinger kernel  
 associated  
to $- \Delta$ in any open subset ${\cal M}$ of $\RR^n$ 
with boundary conditions that  
make the operator positive and self-adjoint, then  
\begin{equation} U_1(t,x,y) = (4 \pi t)^{-n/2} e^{i|x-y|^2/4t}
  + o(t^{\infty}) ,  
 \quad  
\mbox{in} \ {\cal D}'      , \quad \mbox{as} \ t \rightarrow 0^{+}    ,  
\label{(6.20)}   
\end{equation}  
uniformly and strongly on compacts of ${\cal M}$.  In particular,  
\begin{equation} U_1(t,x,x) = (4 \pi t)^{-n/2} + o(t^{\infty})      ,
  \quad  
\mbox{in} \ {\cal D}'      , \quad \mbox{as} \ t \rightarrow 0^{+}    ,  
\label{(6.21)}   
\end{equation}  
uniformly and strongly on compacts of~${\cal M}$.  
 (We reiterate that this holds in the {\em interior\/} of ${\cal M}$ 
only.)  
 $\;\;\bull$  
  \smallskip

Apart from (\ref{(6.21)}) we have had little to say in this section 
about the much-studied {\em diagonal\/} expansions of the heat and
Schr\"odinger kernels.
See, however, Ref.~\onlinecite{EGV}, where the methods of this
paper are applied in that arena.

\section{Expansion of Green kernels II{}: Global expansions}  
\label{sec7}

When considering a second-order differential operator  
 on a one-dimensional manifold,  
the variable $\lambda$ of the spectral density $e(x,y;\lambda)$  
 is often replaced  
by the variable $\omega$ defined by $\omega^2=\lambda$.
    For instance, the   
asymptotic behavior of $\omega_n=\lambda_n^{1/2}$, 
 the square root of the $n$th  
eigenvalue, has a more convenient form than that for $\lambda_n$ 
in the case  
of regular Sturm--Liouville equations.  
  But, does this change of variable have any effect  
on the expansion of the associated Green kernels?  
  
Consider, for instance, the behavior of the kernel  
\begin{equation} e(x,y;\lambda) = \frac{2}{\pi} \sum_{n=1}^{\infty}
  \sin nx\,\sin ny\,
\delta(\lambda-n^2)  ,    
\label{(7.1)}   
\end{equation}  
which we studied before.
   Let $\omega^2=\lambda$, so that  
\begin{equation} \tilde{e}(x,y;\omega) = e(x,y;\omega^2) =
  \frac{1}{\pi} \sum_{n=1}^{\infty}   
\sin nx\, \sin ny\, \delta(\omega-n)  .  
\label{(7.2)}   
\end{equation}  
The behavior of $\tilde{e}(x,y;\omega)$ at infinity can be obtained 
 by studying  the parametric behavior,  
 i.e., $\tilde{e}(x,y;\sigma \omega)$ as $\sigma \rightarrow  
\infty$. 
   Letting $\varepsilon=1/\sigma$, we are led to consider the 
development  of  
\begin{equation} \frac{1}{\pi} \sum_{n=1}^{\infty} 
 \sin nx\, \sin ny\, \phi(\varepsilon n) =   
\left\langle \tilde{e}(x,y;\omega),\,\phi(\varepsilon n) 
 \right\rangle_{\omega} 
\label{(7.3)}   
\end{equation}  
as $\varepsilon \rightarrow 0^{+}$.  
  The even moments of $\tilde{e}$ coincide   
with those  of $e$ of half the order; i.e., 
\begin{equation} \langle \tilde{e}(x,y;\omega),\omega^{2k} \rangle = 
 \delta^{(2k)}(x-y)  .
\label{(7.4)}   
\end{equation}  
But we also have to consider odd-order moments, such as  
\begin{equation} J(x,y) = \langle \tilde{e}(x,y;\omega),\omega \rangle
  = \sum_{n=1}^{\infty}  
n\, \sin nx\,  \sin ny  .  
\label{(7.5)}   
\end{equation}  
The operator corresponding to $J$ is not the derivative $d/dx$,  
 but rather $(d/dx)Q$,  
where $Q$ is a Hilbert transform.  
  Thus, $J$ is not a local operator: in general  
$\mathop{\rm supp}{J(\phi)}$ is not contained in $\mathop{\rm supp}\phi$. 
   Thus, the expansion  
\begin{equation} \tilde{e}(x,y;\sigma\omega) \sim 
 \frac{\delta(x-y)\delta(\omega)}  
{\sigma} - \frac{J(x,y)\delta'(\omega)}{\sigma^2} + \frac{\delta''(x-y)  
\delta''(\omega)}{2\sigma^3} -\cdots   
\label{(7.6)}   
\end{equation}  
 has a nonlocal character:  
 The expansion of $\langle \tilde{e}(x,y;\sigma\omega),  
f(y)\rangle$ as $\sigma \rightarrow \infty$ 
  may have nonzero contributions outside  
of $\mathop{\rm supp} f$.  
  
The change of variable $\omega^2=\lambda$, which seems so innocent,  
 has introduced a  
new phenomenon into the expansion of the spectral density: 
  the appearance of nonlocal terms.  
  This will also apply to the expansion of the corresponding Green   
kernels.  
 Consider the expansion of the cylinder kernel (\ref{(A.9b)})
\begin{equation} T(t,x,y) = \frac{t}{\pi((x-y)^2+t^2)} \, .  
\label{(7.7)}   
\end{equation}  
We have  
\begin{equation} \frac{t}{\pi((x-y)^2+t^2)} \sim
  \delta(x-y) + \frac{t}{\pi(x-y)^2} -   
\frac{\delta''(x-y)t^2}{2} + \frac{t^3}{\pi(x-y)^3} + \cdots \quad
\mbox{as} \ t \rightarrow 0^{+}  ,
\label{(7.8)}   
\end{equation}  
which is a nonlocal expansion because the odd-order terms are nonlocal. 
  In   particular, the expansion does not vanish when $x \neq y$.  
  
Why?  Why do the results of the previous section fail?  
 The operator corresponding  
to (\ref{(7.7)}) is $e^{-tH^{1/2}}$, where $H=-d^2/dx^2$ on $\RR$.   
 Thus, we arrive at the  
same question: How does the change 
 $\lambda\mapsto\lambda^{1/2}=\omega $  affect our results?   
   However, the answer is now clearer:
$ H^{1/2}$ is a {\em nonlocal\/} pseudodifferential operator, and
 $J(x,y)$ is simply its kernel.
 The nonlocality of the odd-order moments just reflects the fact 
that the basic operator $ H^{1/2}$ is not local.
   
The behavior of summability after changes of variables was already 
studied at  the beginning of the century by Hardy\cite{Ha1,HR} 
 and is the central 
 theme in Ref.~\onlinecite{F97}. 
  
When using the distributional approach, we can see that the change 
$\omega=  \lambda^{1/2}$, and similar ones, do not introduce problems 
at infinity.  Rather,  the point is that the change introduces a 
new structure at the origin.  
  
Let $f \in {\cal D}'$ have $\mathop{\rm supp} f \subseteq (0,\infty)$.
    Then $f(\omega^2)$ is a  
well-defined distribution, given by  
\begin{equation} \langle f(\omega^2),\phi(\omega) \rangle = 
 {\textstyle \frac12} \langle f(\lambda),\,\lambda^{-1/2}  
\phi(\lambda^{1/2}) \rangle .  
\label{(7.9)}   
\end{equation}  
When $0 \in\mathop{\rm supp} f$, however, 
 there is no canonical way to define $f(\omega^2)$. 
That there are no problems at infinity follows from the results 
 of  Ref.~\onlinecite{E96}.  
  
\smallskip{\bf Lemma 7.1.} 
{\em Let $f$ be in ${\cal D}'(\RR)$ with 
 $\mathop{\rm supp} f \subseteq (0,\infty)$. 
   Then  
$f$ is distributionally small at infinity if and only if $f(\omega^2)$ is.
 }
 \smallskip 
    
{\em Proof:\/}
  The generalized function $f$ is distributionally small at infinity  
if and only if it belongs to ${\cal K}'$.   
 Thus it suffices to see that $f(\lambda)$  
belongs to ${\cal K}'$ if and only if $f(\omega^2)$ does, and, by duality, 
 it suffices  
to see that if $\,\mathop{\rm supp}\phi \subseteq (0,\infty)$ then 
 $\phi(\omega)$ belongs to  
${\cal K}$ if and only if $\lambda^{-1/2}\phi(\lambda^{1/2})$ does.
  $\;\; \bull$  
\smallskip

Therefore, 
 if $e(x,y;\lambda)$ is the spectral density kernel corresponding to a  
positive self-adjoint operator~$H$,
 then $e(x,y;\omega^2)$ is also distributionally  
small as a function of $\omega$, both distributionally and 
 pointwise on $x \neq y$.   
However, as we have seen, 
 the corresponding moment expansion for $e(x,y;\omega^2)$ will contain  
extra terms.
 These arise as a special case of a general theorem that extends 
the conclusions of Sections \ref{sec5} and~\ref{sec6}
  to the situation where the 
function or distribution $g$ does not have a Taylor expansion at 
the origin.
   
The spaces ${\cal A}\{\phi_n\}$ associated to an asymptotic sequence
  are discussed in Refs.~\onlinecite{EK90,EK94}.
  In particular, if $\alpha_n$ is a sequence with   
$\Re e\; \alpha_n \nearrow \infty$,
  then the space ${\cal K}\{x^{\alpha_n}\}$ consists  
of those smooth functions $\phi$ defined on $(0,\infty)$ 
 that have the behavior  
of the space ${\cal K}$ at infinity but at the origin 
  can be developed in a strong  
expansion  
\begin{equation} \phi(x) \sim a_1 x^{\alpha_1} + a_2 x^{\alpha_2} + 
 a_3 x^{\alpha_3} 
 + \cdots \quad
\mbox{as} \ x \rightarrow 0^{+}  .  
\label{(7.10)}   
\end{equation} 
 The point is that the $\alpha_n$ need not be nonnegative integers. 
  
The functionals $\delta_j \in {\cal K}'\{x^{\alpha_n}\}$ given by  
\begin{equation} \langle \delta_j(x),\phi(x) \rangle = a_j  
\label{(7.11)}   
\end{equation}  
play the role of the traditional delta functions.  
Each $f \in {\cal K}'\{x^{\alpha_n}\}$ 
 admits a generalized moment asymptotic expansion,  
\begin{equation} f(\lambda x) \sim \sum_{j=1}^{\infty} 
 \frac{\mu(\alpha_j)\delta_j(x)}  
{\lambda^{\alpha_j+1}} \quad \mbox{as}
  \ \lambda \rightarrow \infty  ,  
\label{(7.12)}   
\end{equation}  
where $\mu(\alpha_j)=\langle f(x),x^{\alpha_j} \rangle$ are the moments.  
 Therefore,  
if $g \in {\cal K} \{x^{\alpha_n}\}$ then the expansion of 
 $G(t)=\langle f(x),g(tx) \rangle$  
can be obtained from (\ref{(7.12)}) as  
\begin{equation} G(t) \sim \sum_{j=1}^{\infty} \mu(\alpha_j)\,a_j 
t^{\alpha_j}    \quad   
\mbox{as} \ t \rightarrow 0^{+}  ,  
\label{(7.13)}   
\end{equation}  
where $a_j=\langle \delta_j(x),g(x) \rangle$.  
  But following the  ideas of   
Section \ref{sec5} we can define $G(t)=\langle f(x),g(tx) \rangle$ 
 when $g$ is a distribution  
of ${\cal S}'$ with $\mathop{\rm supp} g \subseteq [0,\infty)$, 
 whose behavior at the origin is of  
the form  
 \begin{mathletters} \label{(7.14)}
\begin{equation} g(\varepsilon x) \sim \sum_{j=1}^{\infty} a_j\,
  \varepsilon^{\alpha_j}\, x_{+}^{\alpha_j} \quad   
\mbox{as} \ \varepsilon \rightarrow 0^{+}  ,  
\label{(7.14a)}   
\end{equation}  
a fact that we express by saying that 
 \begin{equation}
 g(x) \sim \sum_{j=1}^{\infty} a_j\, x^{\alpha_j}\,, \quad
 \mbox{distributionally,\quad as}\ x \rightarrow 0^{+}.  
 \label{(7.14b)}\end{equation}   
 \end{mathletters}
  Then $G(t)$ will have the same expansion  
(\ref{(7.13)}), but in the average or distributional sense. 
  A corresponding result for   
operators also holds.  
  
\smallskip
 {\bf Theorem 7.1.} 
{\em Let $H$ be a positive self-adjoint operator on the domain  
${\cal X}$ of the Hilbert space ${\cal H}$.   
 Let ${\cal X}_\infty$ be the intersection of the domains   
of $H^n$ for $n \in \NN$.   
 Let $g \in {\cal S}'$ with $\mathop{\rm supp} g \subseteq [0,\infty)$   
 have a distributional expansion of the type $g(x) \sim   
\sum_{j=1}^{\infty} a_j\, x^{\alpha_j}$ as $x \rightarrow 0^{+}$, where   
$\Re e\, \alpha_n \nearrow \infty$. 
   Then $G(t)=g(tH)$ can be defined as an   
element of ${\cal S}'(\RR,L({\cal X}_{\infty},{\cal H}))$ 
 with support contained in $[0,\infty)$,  
and $G(t)$ admits the distributional expansion  
\begin{equation} G(t) \sim \sum_{j=1}^{\infty} a_j\, H^{\alpha_j} 
\, t^{\alpha_j} \, , \quad   
\mbox{as} \ t \rightarrow 0^{+}  , \quad \mbox{in} \ {\cal D}'  .  
\label{(7.15)}   
\end{equation}  
When $g$ belongs to ${\cal K}\{x^{\alpha_n}\}$, (\ref{(7.15)}) 
 becomes a pointwise expansion.  
$\;\; \bull$
}
\smallskip

Therefore, we may generalize our previous discussion as follows: 
{\em If $H$ is a differential  
operator, then the expansion of the Green function of $g(tH)$ 
is local or global
depending on whether the expansion (\ref{(7.14)})
  of $g$ at the origin is of the 
 Taylor-series type or not.}  
  
\smallskip
 {\bf Example.} The small-$t$ expansion of the cylinder function $T(t,x,y)$   
described in (\ref{(A.3)}) is given by   
\begin{equation}  
T(t,x,y)\sim \sum_{n=0}^\infty {(-1)^n\,
 H^{n/2}_x(\delta (x-y))\,t^n\over n!}
 \quad   {\rm as}\ t\to 0^+.  
\label{(7.16)}  
\end{equation}  
The expansion is pointwise in $t$ and distributional in $(x,y)$. The expansion  
is also pointwise in $t$ for $x\not= y$, 
 but we do not get $T(t,x,y)=o(t^\infty)$  
because the odd terms in the expansion do not vanish for $x\not= y$. 
This type  
of behavior is typical of  harmonic functions near boundaries. 
$\;\;\bull$  
\smallskip

We may look at the locality problem from a different perspective.
  Suppose  
$H_1$ and $H_2$ are two different self-adjoint extensions of the same  
differential operator on a subset $U$ of ${\cal M}$.
  Then the two cylinder kernels  
$T_1(t,x,y)$ and $T_2(t,x,y)$ have different expansions as $t\to 0^+$ even if  
$(x,y)\in U\times U$. 
 The same is true of the associated Wightman functions.  
In both cases the small-$t$ expansion of the Green kernel reflects some global  
properties of the operators $H_1$ and~ $H_2\,$. 
 Since the cylinder and Wightman  
functions are constructed from the operators $H_1^{1/2}$ and $H_2^{1/2}$,
  one  
may ask if this nonlocal character can already be observed in the spectral  
densities $e_{H_j^{1/2}}(x,y;\lambda)$, $j=1,2$. 
 Interestingly, the nonlocal  
character cannot be seen in the Ces\`aro behavior, since according to  
Theorem~7.2 below we have  
\begin{equation}  
e_{H_1^{1/2}}(x,y; \lambda)=e_{H_2^{1/2}}(x,y;\lambda)
 +o(\lambda^{-\infty})\quad (C)  
\quad {\rm as}\ \lambda\to\infty\ ,  
\label{(7.17)}  
\end{equation}  
for $(x,y)\in U\times U$. 
  Instead, the nonlocal character of the small-$t$  
expansion of the Green kernels is explained by the difference in the moments.  
(Recall  Theorem~2.1
and the formula (\ref{(2.8)}).)  
  
We finish by giving a result that justifies (\ref{(7.17)}) and also has an   
interest of its own.

\smallskip{\bf Theorem 7.2.}
{\em  Let $H_1$ and $H_2$ be two pseudodifferential operators  
acting on the manifold ${\cal M}$, 
 with spectral densities $e_j(x,y;\lambda)$ for  
$j=1,2$. 
 Let $U$ be an open set of ${\cal M}$ and suppose that $H_1-H_2$ is a  
smoothing operator in $U$. Then  
\begin{equation}  
e_1(x,y;\lambda)=e_2(x,y;\lambda)+o(\lambda^{-\infty})\quad (C)\quad
   \mbox{as}\   
\lambda\to\infty,  
\label{(7.18)}  
\end{equation}  
in the topology of the space ${\cal E} (U\times U)$ and,
  in particular, pointwise  
on $(x,y)\in U\times U$.
}
\smallskip  
  
{\em Proof:\/} If $\phi\in{\cal D} (\RR)$, 
 then $\phi(H_1)-\phi(H_2)$ is a smoothing  
operator, thus $\langle e_1(x,y;\lambda)-e_2(x,y;\lambda),\,f(x)g(y)\rangle$  
is a well-defined element of ${\cal D}'(\RR)$ given by   
\begin{equation}  
\langle \langle e_1(x,y;\lambda)-e_2(x,y;\lambda),f(x)g(y)\rangle,\,\phi 
 (\lambda)  
\rangle =\langle (\phi(H_1)-\phi(H_2))f,\,g\rangle.  
\end{equation}  
In general this generalized function is not a distributionally small function  
of $\lambda$, 
 but if $\mathop{\rm supp}{f}\subset U$ and 
 $\,\mathop{\rm supp}{g}\subset U$ then all the  
moments   
\begin{equation}  
\langle\langle e_1(x,y;\lambda)-e_2(x,y;\lambda),f(x)g(y)\rangle 
 ,\,\lambda^n\rangle  
=\langle (H_1^n-H_2^n)f,\,g\rangle  
\end{equation}  
exist because $H_1-H_2$ is smoothing in $U$. 
 Therefore $\langle e_1(x,y;\lambda)-  
e_2(x,y;\lambda),\,f(x)g(y)\rangle$ 
belongs to ${\cal K}'(\RR)$;
 that is, it is a  
distributionally small function. Hence,  
\begin{equation}  
\langle e_1(x,y;\lambda)-e_2(x,y;\lambda),\,f(x)g(y)\rangle =  
o(\lambda^{-\infty})\quad (C)\quad {\rm as}\ \lambda\to\infty  
\label{(7.19)}  
\end{equation}  
for each $f,g\in{\cal E}'(U)$, and (\ref{(7.18)}) follows by duality. 
 $\;\;\bull$   
  
\acknowledgments
We thank I. G. Avramidi for helpful comments.

 \appendix
\section*{The simplest one-dimensional examples}

 Let $H$ be a positive, self-adjoint, second-order linear 
differential operator on scalar functions, on a manifold or region 
${\cal M}$.
 We are concerned with distributions of the type $G(t,x,y)$
 ($t\in\RR$, $x\in{\cal M}$, $y\in{\cal M}$)
 that are integral kernels of parametrized operator-valued 
functions of~$H$.
 In particular:
 
 (1) The {\em heat kernel}, $K(t,x,y)$, represents the operator 
$e^{-tH}$, which solves the heat equation
 \begin{equation}-\,{\partial{\Psi}\over\partial t} = H \Psi, 
 \quad \lim_{t\downarrow0}\Psi(t,x) = f(x), \label{(A.1)}
 \end{equation}
 for $(t,x)\in (0,\infty)\times{\cal M}$, by
 $\Psi(t,x) =  e^{-tH} f(x)$.

 (2) The {\em Schr\"odinger propagator}, $U(t,x,y)$, is the kernel 
of $e^{-itH}$, which solves the Schr\"odinger equation
 \begin{equation}i\,{\partial{\Psi}\over\partial t} = H \Psi, 
 \quad \lim_{t\to0}\Psi(t,x) = f(x), \label{(A.2)}\end{equation}
 for $(t,x)\in \RR\times{\cal M}$, by
 $\Psi(t,x) =  e^{-itH} f(x)$.

 (3) Let $T(t,x,y)$ be the integral kernel of the operator 
$e^{-t{\sqrt H}}$, which solves the elliptic equation
 \begin{equation}-\,{\partial{^2\Psi}\over\partial{t^2}} + H \Psi=0, 
 \quad \lim_{t\downarrow0}\Psi(t,x) = f(x), 
 \quad \lim_{t\to+\infty}\Psi(t,x) = 0, \label{(A.3)}\end{equation}
in the infinite half-cylinder  $(0,\infty)\times{\cal M}$, by
 $\Psi(t,x) =  e^{-t{\sqrt H}} f(x)$.
We shall call this the {\em cylinder kernel\/} of $H$.  It
may also be regarded as the heat kernel of the first-order
 pseudo\-differential operator ${\sqrt H}$.

 (4) The {\em Wightman function}, $W(t,x,y)$, is the kernel of
 $(2{\sqrt H})^{-1} e^{-it{\sqrt H}}$.
 This operator solves the wave equation 
\begin{equation} -\,{\partial{^2\Psi}\over\partial{t^2}} =
  H \Psi \label{(A.4)} \end{equation}
with the nonlocal initial data
 \begin{equation}\lim_{t\downarrow 0} \Psi(t,x) = (2{\sqrt H})^{-1}f(x),
  \quad
\lim_{t\downarrow 0} {\partial{\Psi}\over\partial t} (t,x)
  = -\,{i\over2}\,f(x).
  \label{(A.5)} \end{equation}
 The significance of $W$ is that it is the two-point vacuum 
expectation value of a quantized scalar field satisfying the 
 time-independent, linear field equation (\ref{(A.4)}):
 \begin{displaymath}W(t,x,y) = \langle0| \phi(t,x) \phi(0,y) |0\rangle
  \end{displaymath} 
 (e.g., Ref.~\onlinecite{F89}, Chapters 3--5). 

 These four kernels are rather diverse in their asymptotic
behavior as $t$ approaches~$0$ and also in the convergence 
properties of their spectral expansions.
 Most of the relevant mathematical phenomena that distinguish them 
can be demonstrated already in the simplest case,
\begin{equation}H = 
 -\,{\partial{^2}\over\partial{x^2}} \label{(A.6a)} \end{equation}
 with ${\cal M} $ either $\RR$ or a bounded interval $(0,\pi)$.
 (For other one-dimensional ${\cal M}$ see Ref.~\onlinecite{F97}.) 
 We record here the spectral expansion (Fourier transform or series) 
 of each kernel and also its actual functional value 
 (closed form or image sum).

\subsection*{Case ${\cal M} =\RR$}

\nobreak
 \noindent{\em Heat kernel:} \begin{mathletters}\label{(A.7)}
 \begin{equation}K(t,x,y) =  {1\over 2\pi} \int_{-\infty}^\infty 
 e^{ik(x-y)} e^{-k^2 t} \, dk.  \label{(A.7a)} \end{equation}
 \begin{equation}K(t,x,y) =  (4\pi t)^{-1/2} e^{-(x-y)^2/4t}. 
 \label{(A.7b)} \end{equation}\end{mathletters}

 \noindent{\em Schr\"odinger propagator:}
 \begin{mathletters}\label{(A.8)}
 \begin{equation}U(t,x,y) =  {1\over 2\pi} \int_{-\infty}^\infty 
 e^{ik(x-y)} e^{-ik^2 t} \, dk.  \label{(A.8a)} \end{equation}
 \begin{equation}U(t,x,y) = e^{-i(\mathop{\rm sgn} t)\pi/4}\,
(4\pi t)^{-1/2}
  e^{i(x-y)^2/4t}. \label{(A.8b)} \end{equation} \end{mathletters}

 \noindent{\em Cylinder kernel:}
 \begin{mathletters}\label{(A.9)}
 \begin{equation}T(t,x,y) =  {1\over 2\pi} \int_{-\infty}^\infty 
 e^{ik(x-y)} e^{-|k| t} \, dk.  \label{(A.9a)} \end{equation}
 \begin{equation}T(t,x,y) = {t\over\pi}\, {1\over (x-y)^2 + t^2}\,. 
 \label{(A.9b)} \end{equation}\end{mathletters}

 \noindent{\em Wightman function:}
  \begin{equation}W(t,x,y) =  {1\over 4 \pi} \int_{-\infty}^\infty 
 e^{ik(x-y)} {e^{-i|k| t}\over |k|} \, dk.  \label{(A.10)} \end{equation}
This integral is divergent at $k=0$ and does not make sense even as 
a distribution except on a restricted class of test functions.
 This ``infrared'' problem, which is irrelevant to the main issues 
of the present paper, disappears when one (1) goes to higher 
dimension, (2) adds a positive constant (``mass'') to~$H$, or
 (3) takes one or more derivatives of $W$ with respect to any of 
its variables.
 Therefore, we do not list an integrated form of (\ref{(A.10)}).
 (If we gave one, it would be nonunique and would grow 
logarithmically in $(x-y)^2$, thus being useless for forming an 
image sum for (\ref{(A.14)}).)
For more information about infrared complications in simple model 
quantum field theories, see Ref.~\onlinecite{FR} and 
references therein.

 \subsection*{Case ${\cal M} =(0,\pi)$}

\nobreak
 \noindent{\em Heat kernel:}
 \begin{mathletters}\label{(A.11)}
 \begin{equation}K(t,x,y) =  {2\over \pi} \sum_{k=1}^\infty 
\sin (kx) \sin (ky)  e^{-k^2 t} .  \label{(A.11a)} \end{equation}
 \begin{equation}K(t,x,y) =  (4\pi t)^{-1/2} \sum_{N=-\infty}^\infty
\left[ e^{-(x-y-2N\pi)^2/4t} -  e^{-(x+y-2N\pi)^2/4t} \right]  .
  \label{(A.11b)} \end{equation}
 \end{mathletters}

 \noindent{\em Schr\"odinger propagator:}
  \begin{mathletters}\label{(A.12)}
 \begin{equation}U(t,x,y) =  {2\over \pi} \sum_{k=1}^\infty 
\sin (kx) \sin (ky) e^{-ik^2 t} .  \label{(A.12a)} \end{equation}
 \begin{equation}U(t,x,y) = e^{-i(\mathop{\rm sgn} t)\pi/4}\,
(4\pi t)^{-1/2}
\sum_{N=-\infty}^\infty 
 \left[ e^{i(x-y-2N\pi)^2/4t} - e^{i(x+y-2N\pi)^2/4t} \right].
  \label{(A.12b)} \end{equation}
 \end{mathletters}

 \noindent{\em Cylinder kernel:}
  \begin{mathletters}\label{(A.13)}
 \begin{equation}T(t,x,y) =  {2\over \pi} \sum_{k=1}^\infty 
\sin(kx) \sin (ky) e^{-k t}.  \label{(A.13a)} \end{equation}
 \begin{equation}T(t,x,y) = {t\over\pi} \sum_{N=-\infty}^\infty
\left[  {1\over (x-y-2N\pi)^2 + t^2} - {1\over (x+y-2N\pi)^2 + t^2} 
\right].
  \label{(A.13b)} \end{equation}
 In this case a closed form is obtainable:
 \begin{equation}T(t,x,y) = {1\over2\pi} 
 \left[ {\sinh t \over \cosh t - \cos(x-y) }
 - {\sinh t \over \cosh t - \cos(x+y) } \right]. \label{(A.13c)} \end{equation}
  \end{mathletters}

 \noindent{\em Wightman function:}
  \begin{mathletters}\label{(A.14)}
 \begin{equation}W(t,x,y) = -\, {1\over \pi} \sum_{k=1}^\infty 
\sin(kx) \sin(ky)\, {e^{-ik t}\over k}\, . \label{(A.14a)} \end{equation}
 This is expressible in the closed form
 \begin{equation} W(t,x,y) =  {1\over 4\pi} \ln \left| {\cos t - \cos (x+y) 
 \over \cos t - \cos (x-y)} \right|
 + {i\over4} \,P(t,x,y), \label{(A.14b)} 
 \end{equation} \end{mathletters}
 where
 \begin{equation} P(t,x,y) = \cases {
-1 & for $2k\pi -f(x+y) < t < 2k\pi -|x-y|$, \cr
0 & for $2k\pi -|x-y| < t < 2k\pi +|x-y|$, \cr
1& for $2k\pi +|x-y| < t < 2k\pi +f(x+y)$, \cr
0 & for $2k\pi +f(x+y) < t < 2(k+1)\pi -f(x+y)$, \cr} 
\label{(A.15)}\end{equation}
for $k\in \ZZ$, with
\begin{equation}f (z) = \cases{
z &if $ 0\le z\le \pi$, \cr
2\pi - z &if $ \pi \le z \le 2\pi$. \cr}
\label{(A.16)}\end{equation}
($P$ is essentially the standard Green function for the Dirichlet problem
for the wave equation --- the d'Alembert solution modified by reflections.)

\smallskip
{\em Proof of (\ref{(A.13c)}):\/}
Write (\ref{(A.13a)}) as
\begin{displaymath}-\,{1\over 2\pi} \sum_{k=1}^\infty
\Bigl[ e^{ik(x+y)-kt} + e^{-ik(x+y)-kt} - e^{ik(x-y)-kt} -e^{-ik(x-y)-kt}
\Bigr].\end{displaymath}
Evaluate each sum by the geometric series
\begin{equation} \sum_{k=1}^\infty e^{kz} = {e^z \over 1-e^z}\,. 
 \label{(A.17)}\end{equation}
Some algebraic reduction yields (\ref{(A.13c)}).
$\;\;\bull$

\smallskip
{\em Proof of (\ref{(A.14b)}):\/}
Start with the well-known dispersion relation
\begin{equation}{1\over x\pm i0} = 
 {\cal P} {1\over x} \mp \pi i \delta(x), \label{(A.18)}\end{equation}
where $ {\cal P} {1\over x}$ is a principal-value distribution, whose 
antiderivative is $\,\ln |x| + \mbox{constant}$.
This generalizes easily to:

 \smallskip
{\bf Lemma A.1.}
 {\em   Let $F$ be analytic in a region 
$\Omega \setminus \{x_1, \ldots, x_n \}$, $x_j \in \RR$, where $\Omega$
intersects $\RR$ on $(a,b)$ and each $x_j$ is a simple pole of $F$ with
residue $\alpha_j\,$.
Then
\begin{equation} F(x-i0) = {\cal P} F(x) + \pi i
  \sum_{j=1}^n \alpha_j \delta(x-x_j). 
\label{(A.19)}\end{equation}
$\bull$
 }
 \smallskip

To apply this to (\ref{(A.14b)}), replace $t$ by $ib$ in 
 (\ref{(A.13a)}) and (\ref{(A.13c)}):
\begin{displaymath} \sum_{k=1}^\infty \sin(kx) \sin(ky) e^{-ikb} =
{i\over 4} \left( {\sin b \over \cos b - \cos(x-y)}
-  {\sin b \over \cos b - \cos(x+y)} \right),\end{displaymath}
for $\,\mbox{Im\,} b < 0$.
There are poles at $b=2k\pi \pm (x-y)$ with residue $i/4$
and at  $b=2k\pi \pm (x+y)$ with residue $-i/4$.
Thus for $b\in \RR$,
\begin{displaymath} \sum_{k=1}^\infty \sin {kx} \sin {ky}\, e^{-ikb} =
{i\over 4} \left[ {\cal P} \left( {\sin b \over \cos b -\cos(x-y) }
- {\sin b \over \cos b -\cos(x+y) }  \right)\right.\end{displaymath}
\begin{equation}\left. {}+ \pi i \sum_{k=-\infty}^\infty
 \bigl(\delta(b-2k\pi -x+y) +\delta(b-2k\pi +x-y) -
\delta(b-2k\pi -x-y) -\delta(b-2k\pi +x+y) \bigr)
\vphantom{\sin b \over \cos b -\cos(x-y) } \right].
\label{(A.20)}\end{equation}
Also, integrating (\ref{(A.13a)})--(\ref{(A.13c)})
  and letting $t\to 0^+$, one gets
\begin{equation}\sum_{k=1}^\infty {\sin{kx} \sin{ky} \over k} =  {1\over 4}
\ln \left( {1-\cos(x+y) \over 1-\cos(x-y)}\right).
\label{(A.21)}\end{equation}
Therefore, integrating (\ref{(A.20)}) yields
\begin{displaymath}\sum_{k=1}^\infty {\sin{kx} \sin{ky} \,e^{-ikt} \over k} =
-i \int_0^t \sum_{k=1}^\infty \sin{kx} \sin{ky} \,e^{-ikb} \, db
+ \sum_{k=1}^\infty {\sin{kx} \sin{ky}  \over k} \end{displaymath}
\begin{displaymath}{}= {1\over 4} \int_0^t {\cal P}
  \left( {\sin b \over \cos b -\cos(x-y) }
- {\sin b \over \cos b -\cos(x+y) }  \right)\,db
+ {1\over 4}\ln \left( {1-\cos(x+y) \over 1-\cos(x-
y)}\right)\end{displaymath}
\begin{displaymath} {}+ {\pi i\over4} \int_0^t  \sum_{k=-\infty}^\infty
 \bigl(\delta(b-2k\pi -x+y) +\delta(b-2k\pi +x-y) -
\delta(b-2k\pi -x-y) -\delta(b-2k\pi +x+y) \bigr) 
\,db\end{displaymath}
\begin{displaymath}{}= {1\over 4} \ln \left| {\cos t-\cos(x+y) \over 
 \cos t-\cos(x-y)}\right|
+ {\pi i\over 4} \,P(t,x,y),\end{displaymath}
which is (\ref{(A.14b)}).
$\;\;\bull$
 
\end{document}